\documentclass[a4paper,11pt]{article}

\usepackage{amssymb}
\usepackage{amsfonts,dsfont}
\usepackage{amsmath}
\usepackage{amsthm}
\usepackage{cite}
 \usepackage{enumitem}
\numberwithin{equation}{section}

\usepackage{xcolor}

\newcommand{\R}{\mathbb{R}}
\newcommand{\Z}{\mathbb{Z}}
\newcommand{\N}{\mathbb{N}}

\newcommand{\cuad}{{\sqcap\kern-.68em\sqcup}}

\newtheorem{definition}{Definition}[section]
\newtheorem{theorem}[definition]{Theorem}
\newtheorem{proposition}[definition]{Proposition}
\newtheorem{weylproposition}{Proposition}[section]

\newtheorem{lemma}[definition]{Lemma}
\newtheorem{remark}{Remark}[section]
\newtheorem{corollary}{Corollary}[section]

\renewcommand{\phi}{\varphi}

\newcommand{\cH}{{\mathcal H}}

\newcommand{\bH}{{\mathbb H}}

\newcommand{\bX}{{\mathbb X}}

 \allowdisplaybreaks

\begin{document}

\begin{center}
{\bf \large  Dirichlet eigenvalues of Laplacian with potential on general graphs
\\[2mm]
and Weyl limits in integer lattice spaces $\mathbb{Z}^d$
}

\bigskip\bigskip

  {\small  Huyuan Chen\footnote{chenhuyuan@yeah.net, chenhuyuan@simis.cn} \qquad   Lixun Jia\footnote{lxjia26@m.fudan.edu.cn}
 \bigskip

  $ ^{1}$   Center for Mathematics and Interdisciplinary Sciences, Fudan University, \\
Shanghai 200433, China\\[8pt]

   $ ^{1,2}$ Shanghai Institute for Mathematics and Interdisciplinary Sciences,\\
 Shanghai 200433, China\\[12pt]

}
  \medskip

\begin{abstract}
  In this paper, we first establish the existence of weighted Dirichlet eigenvalues of the Laplacian with a vertex weight on an infinite, connected and locally finite graph,  via building compact embeddings of suitable spaces,  when the weight decays faster at infinity than the optimal Hardy weight. In the particular integer lattice space $\mathbb{Z}^d$, we apply a segmentation technique to estimate the mean Riesz potential over the entire space, thereby deriving estimates for the eigenvalue counting function and establishing the Weyl limit of the eigenvalues.

 \end{abstract}
 \end{center}
 {\small
    %\tableofcontents }\vspace{5mm}
   \smallskip

  \smallskip
}
\bigskip

\setcounter{equation}{0}
\section{Introduction}
Let $G=(V,E,\mu)$ be {\it an infinite, connected, locally finite   transient graph,} where $V$ denotes the vertex set, $E$ denotes the set of edges, and $\mu_{xy}:V\times V\to[0,\infty)$ is a symmetric function satisfying $\mu_{xy}=\mu_{yx}$. In particular, $\mu_{xx}=0$ for every $x\in V$. The edge weights induce a vertex measure defined by
$$
\mu(x)=\sum_{y\sim x}\mu_{xy}.
$$
Since the graph is infinite, connected, and locally finite, we have $0<\mu(x)<\infty$ for every $x\in V$.  The  Laplacian is defined by
$$
-\Delta u(x)=\frac{1}{\mu(x)}\sum_{y\sim x}\mu_{xy}\bigl(u(x)-u(y)\bigr),\qquad x\in V.
$$
 Let   $o\in V$ and $G_o$ be the positive minimal Green function
$$-\Delta G_o=\delta_o\quad {\rm in}\ \, V,$$
where $\delta_o$ is the Dirac mass at $o$.  Assume that {\it $G_o$ is proper. }

This paper  concerns  the Dirichlet eigenvalues  of the Laplacian with a weight  on an infinite, connected and locally finite graph $G$,
\begin{equation}\label{E1}
-\Delta u =\lambda W u \quad\  {\rm in}\  \,  V, 
\end{equation}
where   $W:V\to[0,\infty)$ is a nonzero vertex weight.

The study of partial differential equations on graphs reveals profound connections between discrete geometry and the analytic properties of solutions. In recent years, this area has attracted considerable attention, see  \cite{CHZ,CH,CH1,CH2,LLY2011, GLY2016, GLY2017, HM2015, LW2017, KL2012}.
Eigenvalue problems on graphs have been extensively investigated, see \cite{BGK2015,GHKLW2015,KL2010,CGY2000,LOGT2014}.
In particular, eigenvalue problems for graph Laplacians serve as a bridge linking discrete geometry with diffusion processes and spectral analysis \cite{KL2010}. In the discrete spectrum setting, eigenvalues govern the decay rates of associated diffusion modes, while their distribution encodes structural information such as the underlying graph topology and edge weights \cite{KL2010,BGK2015}. Spatial inhomogeneity can be incorporated either via an additive potential, yielding a discrete Schr\"odinger operator \cite{BGK2015}, or through positive vertex weights in a generalized eigenvalue problem which is equivalent to modifying the underlying vertex measure \cite{GHKLW2015}.
 Chung \cite{CGY2000} developed systematic methods that connect geometric conditions such as isoperimetric and Sobolev inequalities to eigenvalue estimates, extending classical results from manifolds to weighted graphs. Complementing this approach, Lee \cite{LOGT2014} established higher order Cheeger type inequalities for finite graphs, relating small eigenvalues of the normalized Laplacian to the existence of multiple nearly disconnected subgraphs.

When the  graph is infinite, the first fundamental question in spectral theory concerns the choice of a self adjoint realization of the Laplacian. Haeseler \cite{HKLW2012} systematically investigated Dirichlet and Neumann realizations of graph Laplacians, elucidating their relationship within the framework of nonnegative self adjoint operators and Dirichlet forms. Their work establishes a rigorous operator theoretic foundation for spectral analysis on infinite graphs.
The second key issue is that the spectrum need no longer be discrete. Consequently, one must address the problem of establishing compact embeddings between appropriate function spaces. A natural approach to proving compactness is to combine structural properties of the energy space such as coercivity or Poincar\'{e} type inequalities with finiteness conditions on the vertex measure.

Moreover,  Hua et al. \cite{HKSW2023} derived lower bounds on eigenvalue growth by leveraging Sobolev type inequalities on graphs endowed with finite measure. These results highlight how functional inequalities and measure theoretic constraints shape the spectral properties of graph Laplacians.   Our purpose in this article is to consider the properties of
weighted Dirichlet eigenvalues for the Laplacian on  an infinite, connected, locally finite   transient graph.

Before stating our main results,  we first introduce some notations.  Let
$$
C_c(V)
=
\left\{u:V\to\mathbb R:
\operatorname{supp}u\text{ is finite}\right\},
\qquad
\operatorname{supp}u=\{x\in V:u(x)\neq0\},
$$
and
$$
\ell_W^2=\ell^2(V,W\mu)
=
\Big\{u:V\to\mathbb R:
\sum_{x\in V}\mu(x)W(x)|u(x)|^2<\infty\Big\}.
$$
Set
$$
(u,v)_{\ell_W^2}=
\sum_{x\in V}\mu(x)W(x)u(x)v(x),$$
then $$
\|u\|_{\ell_W^2}=(u,u)_{\ell_W^2}^{\frac{1}{2}}.
$$
For $u,v\in C_c(V)$, define the energy inner product and its associated norm by
$$
(u,v)_{1}
=
\frac{1}{2}\sum_{x,y\in V}
\mu_{xy}\bigl(u(x)-u(y)\bigr)\bigl(v(x)-v(y)\bigr),
\qquad
\|u\|_{1}=(u,u)_{1}^{\frac{1}{2}}.
$$
The inner products $(\cdot,\cdot)_{1}$ and $(\cdot,\cdot)_{\ell_W^2}$ extend uniquely by density to the completion 
$$\cH^1:=\overline{C_c(V)}^{\,\|\cdot\|_{1}}.$$  {\it A pair $(\lambda,u)\in \R\times\cH^1$ is said to be a solution of (\ref{E1}) if
$$
(u,\phi)_{1}=\lambda (u,\phi) _{\ell_W^2}
\quad\ {\rm for\ all}\  \phi\in\cH^1.
$$}

 Note that the optimal Hardy weight on $V$ is investigated  in  \cite[Theorem 0.2]{KellerPinchoverPogorzelski}. Precisely,
the Hardy inequality writes that for some  $C>0$
\begin{equation}\label{Hardy-inequality}
\sum_{x\in V}\mu(x)H_o(x)|u(x)|^2
\leq C\|u\|_1^2,
\qquad \forall\, u\in C_c(V),
\end{equation}
where
\begin{equation}\label{Hdef}
H_o(x)=\frac{-\Delta (\sqrt{G_o(x)})}{\sqrt{G_o(x)}}\quad {\rm for\ all}\ \, x\in V.
\end{equation}

  Let $W:V\to[0,\infty)$ be a nonzero function satisfying one of the following conditions.
\\
\\
$(\mathbf{W_0})$ There exists a constant $M>0$ such that $W(x)\le M H_o(x)$, $\forall x\in V$.\\
\\
$(\mathbf{W_1})$   $\frac{W(x)}{H_o(x)}\to 0$ as $d_G(o,x)\to+\infty$, where $o\in V$ and  the graph distance is defined by
$$
d_G(x,y)
=
\min\big\{m: \exists x_0,\ldots,x_m\in V\text{ with }x_0=x,\ x_m=y, x_{j-1}\sim x_j\text{ for }j=1,\ldots,m \big\}.
$$
Under above assumptions on the weight $W$, we have the following embedding results.

\begin{proposition}
Assume that  $W:V\to[0,\infty)$ is a nonzero vertex weight.

\textup{(i)} If $(\mathbf{W_0})$ holds, then the embedding of $\cH^1\hookrightarrow\ell_W^2$ is continuous. In particular, there exists a constant $C>0$ such that
$$
\sum_{x\in V}\mu(x)W(x)|u(x)|^2\leq C\|u\|_1^2,\quad u\in\cH^1.
$$

\textup{(ii)} If $(\mathbf{W_1})$ holds, then the embedding of $\cH^1\hookrightarrow\ell_W^2$ is compact. In particular, the same embedding inequality holds for some $C>0$.

\textup{(iii)} If
$$\sup_{\substack{u\in C_c(V)\\u\neq0}}\frac{\displaystyle\|u\|^2_{\ell_W^2}}{\|u\|_1^2}=+\infty,$$
then $\cH^1\not\subset \ell_W^2$.
\end{proposition}

Under the assumption $(\mathbf{W_1})$,  we can obtain the existence of the following eigenvalue sequence.

\begin{theorem}\label{Teo 1}
Assume that    $W:V\to[0,\infty)$ is a nonzero vertex weight  satisfying $(\mathbf{W_1})$, then problem (\ref{E1})
 admits  a sequence of eigenvalues
 $$0<\lambda_1\leq\lambda_2\leq\cdots$$ and corresponding eigenfunctions $\{\xi_k\}_{k\in\mathbb{N}}\subset \cH^1$, which have the following properties:
\\
(i)  $\lambda_1>0$ and
$$
\lambda_k=\min\big\{\|u\|_{1}^2:\ u\in \cH^1_k,\ \|u\|_{\ell_W^2}=1\big\},
$$
where $\cH^1_1=\cH^1$ and
$$
\cH^1_k=\big\{u\in\cH^1:(u,v)_{\ell_W^2}=0\text{ for every }v\in \operatorname{span}\{\xi_1,\ldots,\xi_{k-1}\}\big\}.
$$
(ii) Let $\#\operatorname{supp} W$  be the number of points of $\operatorname{supp} W$.\\
If $\#\operatorname{supp}W=m<\infty$, we have exactly $m$ eigenvalues.\\
 If $\#\operatorname{supp}W=+\infty$, the set of eigenvalues is infinite and
 $$\lambda_k\to+\infty\quad {\rm as} \ \, k\rightarrow\infty. $$
(iii)
$$
(\xi_i,\xi_j)_{\ell_W^2}=\delta_{ij},
\qquad
(\xi_i,\xi_j)_{1}=\lambda_i\delta_{ij}.
$$
Moreover, the eigenfunctions' set $\{\xi_k\}_{k}$  forms an orthonormal basis of $\ell_W^2$.
\end{theorem}

Next we concentrate on the Dirichlet eigenvalues on the integer lattice graph. Let $\Z^d$ be the \textit{d-dimensional integer lattice graph},   which consists of the set of vertices
$$\mathbb{Z}^d=\big\{x=(x_1,\cdots, x_d):x_1,\cdots, x_d\in \Z\big\}$$ and the set of edges $E=\big\{\{x,y\}: x,y\in \Z^d, |x-y|=1\big\},$    together with the edge weight  defined by $\omega:\Z^d\times \Z^d\to [0,+\infty), $
%\[\omega:\Z^d\times \Z^d\to [0,+\infty),\]
\begin{equation*}
    \omega_{xy}= \begin{cases}
1\ \,\, \text{ if } \displaystyle |x-y|:=  \sum^d_{k=1}|x_k-y_k|=1,\\
0\ \,\, \text{ otherwise}
\end{cases}
\end{equation*}
and
$$x\sim y\quad {\rm if}\quad \omega_{xy}=1. $$
The Laplace operator on $\Z^d$  is defined as
 $$\Delta_{\Z^d} u(x)=\sum_{y\sim x}\big(u(y)-u(x)\big)\quad\ \text{ for all }x\in\mathbb{Z}^d.$$

 The optimal Hardy inequality  from \cite[Theorem 7.1]{KellerPinchoverPogorzelski} (also see \cite{KPP1}) has
the form that
\begin{equation}\label{Htineq}
\frac{1}{2}\sum_{x\in\Z^d}  \sum_{y\sim x} \big( \varphi(x)-\varphi(y)\big)^2 \geq \sum_{z\in\Z^d}  \overline{H}_0(z)\varphi(z)^2 \quad{\rm for\ all}\  \varphi \in C_c(\Z^d),
\end{equation}
where
\begin{equation}\label{Green 1}
\overline{H}_0(x) =-\big[\Delta\Phi_{d}^{\frac12}(x)\big]\Phi_{d}^{-\frac12}(x)  >0
\end{equation}
and $\Phi_{d}$ is the fundamental solution of Laplacian in $\Z^d$.
Furthermore,  there holds
\begin{equation}\label{ome-1}
\overline{H}_0(x)=\frac{(d-2)^2}{4}|x|^{-2}+O(|x|^{-3})\quad {\rm as}\ \, |x|\to+\infty,
\end{equation}
Then assumption $(\mathbf{W_1})$ reduces to
\begin{align}\label{cond tilde W}
\lim_{  |x|\to+\infty}  W(x)|x|^2= 0.
\end{align}

 In the integer lattice space, we have the following corollary:

\begin{corollary}\label{cr-1}
Let $d\geq3$,  $W\!:\Z^d\to[0,\infty)$ be a nonzero vertex weight  satisfying (\ref{cond tilde W}).
Then problem
\begin{equation}\label{E2}
-\Delta u =\lambda W u \quad {\rm in}\  \Z^d, 
\end{equation}
 admits  a sequence of eigenvalues
 $$0<\lambda_1\leq\lambda_2\leq\cdots$$ and corresponding eigenfunctions $\{\xi_k\}_{k\in\mathbb{N}}\subset \cH^1$, which have the following properties:
\\
(i)  $\lambda_1>0$ and
$$
\lambda_k=\min\big\{\|u\|_{1}^2:\ u\in \cH^1_k,\ \|u\|_{\ell_W^2}=1\big\},
$$
where $\cH^1_1=\cH^1$ and
$$
\cH^1_k=\big\{u\in\cH^1:(u,v)_{\ell_W^2}=0\text{ for every }v\in \operatorname{span}\{\xi_1,\ldots,\xi_{k-1}\}\big\}.
$$
(ii) Let $\#\operatorname{supp} W$  be the number of points of $\operatorname{supp} W$.\\
If $\#\operatorname{supp}W=m<\infty$, we have exactly $m$ eigenvalues.\\
 If $\#\operatorname{supp}W=+\infty$, the set of eigenvalues is infinite and
 $$\lambda_k\to+\infty\quad {\rm as} \ \, k\rightarrow\infty. $$
(iii)
$$
(\xi_i,\xi_j)_{\ell_W^2}=\delta_{ij},
\qquad
(\xi_i,\xi_j)_{1}=\lambda_i\delta_{ij}.
$$
Moreover, the eigenfunctions' set $\{\xi_k\}_{k}$  forms an orthonormal basis of $\ell_W^2$.
\end{corollary}

Note that when $\#\operatorname{supp}W=+\infty$, the number of the eigenvalues is infinite and an interesting question is {\it "what is  the Weyl's  limit?" }

Let's recall the classical  Weyl's limit for the Dirichlet eigenvalues of the Laplacian on bounded domain of the Euclidean space.   For a bounded regular domain $\Omega\subset\mathbb{R}^d$, the Dirichlet   problem
$$
\begin{cases}
-\Delta u=\Lambda u & \text{in } \Omega,\\
\qquad\!\! u=0 & \text{on } \partial\Omega
\end{cases}
$$
has a sequence of eigenvalues
$$0<\Lambda_1(\Omega)\leq \Lambda_2(\Omega)\leq \cdots\leq  \Lambda_k(\Omega)\leq \cdots \to +\infty.  $$
The classical Weyl law \cite{Ivrii2016} gives
$$
\Lambda_k(\Omega)
\sim (c_d|\Omega|)^{-\frac{2}{d}}k^{\frac{2}{d}}
\quad \text{as } k  \to +\infty,
$$
where
$$c_d=(2\pi)^{-d}|B_1(0)|.  $$
Rozenblum \cite{Rozenblum1972} studied the problem on unbounded domains, if $\Omega$ is unbounded and $0<|\Omega|<\infty$, the problem
$$
-\Delta u=\Lambda u \quad \text{in } \Omega,
\qquad u\in \bH_0^1(\Omega),
$$
has discrete spectrum and satisfies the same leading asymptotic
$$
\Lambda_k(\Omega)
\sim (c_d|\Omega|)^{-\frac{2}{d}} k^{\frac{2}{d}}
\quad \text{as } k \to +\infty.
$$
In recent years,  the Weyl's limit   of eigenvalues has been studied extensively, see \cite{Golenia2014, BG2015, BGK2015, BGK2025, RS2009}.

 Our final purpose is to consider  the Weyl limit of the eigenvalues $\lambda_k$ of (\ref{E2}).
  To this end, we need  following notations. Let
$$
\sigma(\theta)
=
2\sum_{j=1}^{d}(1-\cos\theta_j),
\qquad
\theta\in\mathbb{T}^{d}=(-\pi,\pi]^d.
$$
For $t\geq 0$, define
$$
\nu(t)
=
\frac{1}{(2\pi)^d}
\int_{\mathbb{T}^{d}}
\mathbf{1}_{\{\sigma(\theta)\leq t\}}\,d\theta
$$
and
$$
\rho(t)
=
\frac{1}{(2\pi)^d}
\int_{\mathbb{T}^{d}}
\bigl(t-\sigma(\theta)\bigr)_+\,d\theta,
$$
where $a_+=\max\{a,0\}$.
\begin{remark}
 The function $\rho(t)$ satisfies
$$
\rho(t)\sim
\frac{2\omega_d}{(d+2)(2\pi)^d}\,t^{1+\frac{d}{2}}
\quad\text{as }t\rightarrow0^+,
$$
and
$$
\rho(t)=t-2d
\qquad\text{for every }t\geq4d.
$$
\end{remark}

The Weyl's limit is stated as follows.

\begin{theorem}
Let $d\geq3$, $\alpha>2$, and $\beta\in\mathbb{R}$.
Let $W:\mathbb{Z}^d\to[0,\infty)$ satisfy
$$
\frac{W(x)}{|x|^{-\alpha}(\log |x|)^\beta}
\to 1 \quad\text{as } |x|\to+\infty.
$$
Let $\{\lambda_k\}_{k\in\N}$ be the set of Dirichlet eigenvalues of  (\ref{E2}) derived in Corollary \ref{cr-1},  then
\begin{equation}\label{eq w-1}
\lambda_k
\sim
d^\beta C_{d,\alpha} k^{\frac{\alpha}{d}}
(\log k)^{-\beta} \quad\text{as }k\rightarrow+\infty,
\end{equation}
where
\begin{equation}\label{eq Cdef}
C_{d,\alpha}=\Big((1+\frac{d}{\alpha})\int_{\mathbb{R}^d}|y|^\alpha\rho(|y|^{-\alpha})\,dy\Big)^{-\frac{\alpha}{d}}.
\end{equation}
\end{theorem}

 From the Weyl's limit  (\ref{eq w-1}), the power exponent $\alpha>2$ governs the leading order growth rate of the eigenvalues  and  the logarithmic exponent $\beta\in \R$ also plays an important role.
The Weyl's limit   is derived by  using the spatial segmentation and frozing  weights to transform the problem into a local periodic model,  then by controlling various errors through Riesz mean, the exact asymptotic behavior of the counting function and eigenvalues is ultimately restored.  
Furthermore, although the eigenvalue existence condition (\ref{cond tilde W}) in Corollary \ref{cr-1} suggests a natural extension to the critical case $\alpha = 2$ with $\beta < 0$, such an extension is precluded by the divergence of the dimensional constant $C_{d,\alpha}$ defined in (\ref{eq Cdef}): specifically, $C_{d,\alpha} = +\infty$ when $\alpha = 2$. Consequently, the existence of eigenvalues in this regime remains unresolved and is left as an open problem.

\smallskip

The structure of this paper is as follows. In Section 2, we investigate the properties of eigenvalues of matrix,  embedding results, and related estimates. In Section 3, we establish the existence of eigenvalues  and related eigenfunctions.    Finally, Section 4 focuses on Weyl's limit of the eigenvalues in the lattice graph $\Z^d$ with suitable weight $W$.

\section{Preliminaries }
Let two vertices $x,y\in V$ be adjacent, written $x\sim y$ if $\mu_{xy}>0$ and the graph is connected if any two vertices can be joined by a finite path. For a subset $S\subset V$,  its inner vertex boundary are defined by
$$
\delta S=\{x\in S:\text{ there exists }y\in V\setminus S\text{ such that }x\sim y\},
$$
its outer vertex boundary are defined by
$$
\partial S=\{y\in V\setminus S:\text{ there exists }x\in S\text{ such that }x\sim y\}.
$$

Assume that $r:V\to(0,\infty)$ is proper, meaning that
$$
\big\{x\in V:s\leq r(x)\leq t\big\}
\quad\text{is finite for every }0<s\leq t<\infty,
$$
and
$$
\sup_{\substack{x,y\in V\\x\sim y}}
\frac{r(x)}{r(y)}<+\infty.
$$

Let $\bX$ be either $\cH^1$ or $\ell_W^2$. For a sequence $\{u_n\}\subset \bX$ and $u\in \bX$, we say that $u_n$ converges weakly to $u$ in $\bX$, written $u_n\rightharpoonup u$ in $\bX$, if
$$
(u_n,v)_{\bX}\rightarrow (u,v)_{\bX} \quad\text{for every }v\in \bX \quad\text{as }n\rightarrow\infty.
$$
We say that $u_n$ converges strongly to $u$ in $\bX$, written $u_n\rightarrow u$ in $\bX$, if
$$
\|u_n-u\|_{\bX}\rightarrow 0\quad\text{as }n\rightarrow\infty.
$$

\subsection{Eigenvalues of Matrix}

Let $A\in\mathbb{C}^{m\times m}$ be a self adjoint matrix, and define the Counting function and Riesz mean respectively,
\begin{equation}\label{definitionNR}
\mathcal{N}_A(\Lambda)=\#\{j:\lambda_j(A)\leq\Lambda\},\qquad \mathcal{R}_A(\Lambda)=\sum_{j=1}^{m}(\Lambda-\lambda_j(A))_+,
\end{equation}
where $\lambda_j(A)$ are the eigenvalues of $A$, counted with multiplicity.
\\

The following proposition controls the effect of a finite rank perturbation.

\begin{proposition}\label{prop:weyl-rank}
Let $S,T\in\mathbb{C}^{m\times m}$ be self-adjoint matrices, and suppose that
$$
\operatorname{rank}(S-T)\leq r,
\qquad r\in\{0,\ldots,m\}.
$$
Write their eigenvalues in nondecreasing order,
$$
\lambda_1(S)\leq\cdots\leq\lambda_m(S),
\qquad
\lambda_1(T)\leq\cdots\leq\lambda_m(T).
$$
Then the following statements hold.\\
(i) $\lambda_{j-r}(S)\leq\lambda_j(T)\leq\lambda_{j+r}(S)$, where $1\le j+r\le m$ and $1\le j-r\le m$.\\
(ii) For every $\Lambda\in\mathbb{R}$, $|\mathcal{N}_S(\Lambda)-\mathcal{N}_T(\Lambda)|\leq r$.\\
(iii) If $S$ and $T$ are nonnegative, then for every $\Lambda\geq0$, $|\mathcal{R}_S(\Lambda)-\mathcal{R}_T(\Lambda)|\leq\Lambda r$.
\end{proposition}
Before proving, we recall that
$$
\sigma(\theta)
=
2\sum_{j=1}^{d}(1-\cos\theta_j),
\qquad
\theta\in\mathbb{T}^{d}=(-\pi,\pi]^d.
$$
For $t\geq 0$,
$$
\nu(t)
=
\frac{1}{(2\pi)^d}
\int_{\mathbb{T}^{d}}
\mathbf{1}_{\{\sigma(\theta)\leq t\}}\,d\theta
$$
and
$$
\rho(t)
=
\frac{1}{(2\pi)^d}
\int_{\mathbb{T}^{d}}
\bigl(t-\sigma(\theta)\bigr)_+\,d\theta,
$$
where $a_+=\max\{a,0\}$. And we propose several lemmas.
\begin{lemma}\label{num:symbol}
There exists a constant $C_d>0$ such that
$$
\nu(t)\leq C_d\min\{1,t^{\frac{d}{2}}\},
\qquad
\rho(t)\leq C_d\min\{t,t^{1+\frac{d}{2}}\},
\qquad t\geq 0.
$$
\end{lemma}

\begin{proof}
Since the volume of $\mathbb{T}^{d}$ is $(2\pi)^d$, then we have
\begin{equation}\label{v1}
0\leq \nu(t)\leq 1.
\end{equation}
By the Taylor expansion
$$
\cos s
=
1-\frac{s^2}{2}+\frac{s^4}{24}+O(s^6)\quad\mathrm{as}\  s\to 0,
$$
we obtain
$$
\sigma(\theta)
=
\sum_{j=1}^{d}\theta_j^2+O(|\theta|^4)
=
|\theta|^2+O(|\theta|^4)\quad\mathrm{as}\ \theta\to 0.
$$
Consequently, there exist constants
$\delta>0$ and $c_1,c_2>0$, when $|\theta|<\delta$, then
\begin{equation}\label{c1}
c_1|\theta|^2
\leq \sigma(\theta)
\leq c_2|\theta|^2,
\end{equation}
Moreover, $\sigma(\theta)$ is continuous and $\sigma(\theta)=0$ only at
$\theta=0$ on $\mathbb{T}^{d}$. Hence we can choose a $t_0$ is small enough such that if $0\leq t\leq t_0$ and $\sigma(\theta)\leq t$, we have $|\theta|<\delta$.
Then by \eqref{c1} we  have $c_1|\theta|^2\leq t$, i.e., $|\theta|\leq \sqrt{\frac{t}{c_1}}$. It follows that
$$
\{\theta\in\mathbb{T}^{d}:\sigma(\theta)\leq t\}
\subseteq
B_{\sqrt{\frac{t}{c_1}}}(0).
$$
Writing $\omega_d$ for the volume of the unit ball
in $\mathbb{R}^{d}$, we obtain
\begin{equation}\label{v2}
\nu(t)
\leq
\frac{1}{(2\pi)^d}
\left|B_{\sqrt{\frac{t}{c_1}}}(0)\right|
 =\frac{\omega_d}{(2\pi)^d}
\left(\frac{t}{c_1}\right)^{\frac{d}{2}}
 =K_d t^{\frac{d}{2}}
\quad \rm{as}\ 0\leq t\leq t_0,
\end{equation}
where $K_d=\frac{\omega_d c_{1}^{-\frac{d}{2}}}{(2\pi)^d}$. Set
$$
m=\min\{t_0,1\},
\qquad
C_d=\max\{K_d,m^{-\frac{d}{2}},1\}.
$$
For $0\leq t\leq m$, by \eqref{v2} we have
$$
\nu(t)\leq C_d t^{\frac{d}{2}}.
$$
For $m<t\leq 1$, by \eqref{v1} we have
$$
\nu(t)\leq 1
\leq m^{-\frac{d}{2}}t^{\frac{d}{2}}
\leq C_d t^{\frac{d}{2}}.
$$
For $t\geq 1$, by \eqref{v1} we have $\nu(t)\leq 1$. Therefore,
$$
\nu(t)\leq C_d\min\{1,t^{\frac{d}{2}}\},
\qquad t\geq 0.
$$
Finally, for every $a\geq 0$,
$$
(t-a)_+
=
\int_0^t \mathbf{1}_{\{a\leq s\}}\,ds,
$$
so we have
\begin{equation}\label{r1}
\rho(t)=\int_0^t \nu(s)\,ds.
\end{equation}
Since $\nu$ is nondecreasing,
$$
\rho(t)
\leq t\,\nu(t)
\leq C_d\,t\min\{1,t^{\frac{d}{2}}\}
=
C_d\min\{t,t^{1+\frac{d}{2}}\}.
$$
This completes the proof.
\end{proof}

\begin{lemma}
Let $A\in\mathbb{C}^{m\times m}$ be a self adjoint matrix. Then $A$ has real eigenvalues
$$
a_1\leq a_2\leq\cdots\leq a_m
$$
and an orthonormal basis of eigenvectors
$v_1,\ldots,v_m$ satisfying
$$
Av_j=a_jv_j,
\qquad j=1,\ldots,m.
$$
If
$$
u=\sum_{j=1}^{m}c_jv_j,
\qquad c_j\in\mathbb{C},
$$
then
$$
\|u\|_2^2=\sum_{j=1}^{m}|c_j|^2,
\qquad
u^*Au=\sum_{j=1}^{m}a_j|c_j|^2.
$$
For $u\neq 0$, define the Rayleigh quotient by
$$
Q_A(u)
=
\frac{u^*Au}{\|u\|_2^2}
=
\frac{\displaystyle\sum_{j=1}^{m}a_j|c_j|^2}
     {\displaystyle\sum_{j=1}^{m}|c_j|^2}.
$$
Here $u^*$ denotes the conjugate transpose of $u$.
\end{lemma}

\begin{lemma}\label{num:rank-nullity}
Let $D\in\mathbb{C}^{m\times m}$. Define
$$
\ker D
=
\{u\in\mathbb{C}^{m}:Du=0\},
\qquad
\operatorname{Ran}D
=
\{Du:u\in\mathbb{C}^{m}\}.
$$
Then
$$
\dim\ker D+\dim\operatorname{Ran}D=m.
$$
In particular, if
$$
\operatorname{rank}D
=
\dim\operatorname{Ran}D
\leq r,
$$
then
$$
\dim\ker D\geq m-r.
$$
\end{lemma}

\begin{lemma}\label{num:dimension}
Let $U$ and $K$ be subspaces of an $m$-dimensional
vector space. Then
$$
\dim(U+K)
=
\dim U+\dim K-\dim(U\cap K).
$$
Since $\dim(U+K)\leq m$, it follows that
$$
\dim(U\cap K)
\geq
\dim U+\dim K-m.
$$
\end{lemma}

\begin{lemma}\label{num:minmax}
Let $A\in\mathbb{C}^{m\times m}$ be self-adjoint,
with eigenvalues
$$
a_1\leq a_2\leq\cdots\leq a_m,
$$
then for each
$k=1,\ldots,m$,
$$
a_k
=
\min_{\substack{L\subseteq\mathbb{C}^{m}\\ \dim L=k}}
\ \max_{\substack{u\in L\\u\neq 0}}
Q_A(u).
$$
\end{lemma}
Next we prove Proposition~\ref{prop:weyl-rank}.\\
\\
$\mathbf{The\ proof\ of\ Proposition~\ref{prop:weyl-rank}}$
For the self-adjoint matrices $S,T\in\mathbb{C}^{m\times m}$, set
$$
D=S-T,
\qquad
K=\ker D.
$$
(i)
Fix $j$ with $1\leq j+r\leq m$. Let $v_1^S,\ldots,v_m^S$ be an orthonormal basis of eigenvectors of $S$, ordered according to its eigenvalues, and let
$$
U
=
\operatorname{span}
\{v_1^S,\ldots,v_{j+r}^S\}.
$$
Then $\dim U=j+r$. By Lemma~\ref{num:minmax} we have
\begin{equation}\label{Qs}
Q_S(u)\leq\lambda_{j+r}(S),
\qquad 0\neq u\in U,
\end{equation}
where $Q_S(u)$ is the Rayleigh quotient of $S$ restricted to $U\setminus\{0\}$. By the assumption $\operatorname{rank}(S-T)\le r$, Lemma~\ref{num:rank-nullity} gives
\begin{equation}\label{dimk}
\dim K
=
m-\operatorname{rank}(S-T)
\geq m-r.
\end{equation}
Since $\dim U=j+r$ and \eqref{dimk}, the dimension estimate in Lemma~\ref{num:dimension} yields
$$
\dim(U\cap K)
\geq \dim U+\dim K-m
\geq (j+r)+(m-r)-m
=j.
$$
Thus we can choose a $j$-dimensional subspace $Y\subseteq U\cap K$. For every $0\neq u\in Y$, since $Y\subseteq K$ we have
$$
Su=Tu,
$$
and hence
$$
Q_T(u)=Q_S(u)
$$
restricted to $Y\setminus\{0\}$. Besides, since $Y\subseteq U$ and \eqref{Qs}, we obtain
$$
Q_T(u)=Q_S(u)\leq\lambda_{j+r}(S),
\qquad 0\neq u\in Y.
$$
It follows from Lemma \ref{num:minmax} that
\begin{equation}\label{x2}
\begin{aligned}
\lambda_j(T)
&=
\min_{\substack{F\subseteq\mathbb{C}^m\\
                 \dim F=j}}
\ \max_{\substack{u\in F\\u\neq0}}Q_T(u)\\
&\leq
\max_{\substack{u\in Y\\u\neq0}}Q_T(u)\\
&\leq\lambda_{j+r}(S).
\end{aligned}
\end{equation}
Interchanging $S$ and $T$, we similarly obtain
\begin{equation}\label{x3}
\lambda_k(S)\leq\lambda_{k+r}(T),
\qquad 1\leq k+r\leq m.
\end{equation}
Taking $k=j-r$, it follows from \eqref{x2} and \eqref{x3}, we have
\begin{equation}\label{lamuda}
\lambda_{j-r}(S)\leq\lambda_j(T)\le \lambda_{j+r}(S).
\end{equation}
(ii) Fix $\Lambda\in\mathbb{R}$, and set $p=\mathcal{N}_S(\Lambda)=\#\{j:\lambda_j(S)\leq\Lambda\}$. If $p>r$, by the definition of $\mathcal{N}_S(\Lambda)$ we have
$$
\lambda_p(S)\leq\Lambda.
$$
Applying the upper bound proved in part (i) or \eqref{lamuda} with $j=p-r$, we obtain
$$
\lambda_{p-r}(T)\leq\lambda_p(S)\leq\Lambda.
$$
Consequently, by the definition of $\mathcal{N}_T(\Lambda)=\#\{j:\lambda_j(T)\leq\Lambda\}$,
$$
\mathcal{N}_T(\Lambda)\geq p-r=\mathcal{N}_S(\Lambda)-r,
$$
so
$$
\mathcal{N}_S(\Lambda)-\mathcal{N}_T(\Lambda)\leq r.
$$
If $p\leq r$, the same inequality follows directly from
$$
\mathcal{N}_S(\Lambda)-\mathcal{N}_T(\Lambda)\leq \mathcal{N}_S(\Lambda)=p\leq r.
$$
Interchanging $S$ and $T$ gives
$$
\mathcal{N}_T(\Lambda)-\mathcal{N}_S(\Lambda)\leq r.
$$
Therefore,
\begin{equation}\label{NS-NT}
|\mathcal{N}_S(\Lambda)-\mathcal{N}_T(\Lambda)|\leq r.
\end{equation}
\\
(iii) Assume now that $S,T\geq0$, and let $\Lambda\geq0$. For every $a\geq0$,
$$
(\Lambda-a)_+=\int_0^\Lambda\mathbf{1}_{\{a\leq t\}}\,dt.
$$
Since the sums are finite, we may interchange summation and integration. Thus
$$
\begin{aligned}
\mathcal{R}_S(\Lambda)&=\sum_{j=1}^{m}(\Lambda-\lambda_j(S))_+\\
&=\int_0^\Lambda\sum_{j=1}^{m}\mathbf{1}_{\{\lambda_j(S)\leq t\}}\,dt\\
&=\int_0^\Lambda \mathcal{N}_S(t)\,dt.
\end{aligned}
$$
Similarly,
$$
\mathcal{R}_T(\Lambda)=\int_0^\Lambda \mathcal{N}_T(t)\,dt.
$$
Using part (ii) or \eqref{NS-NT}, we conclude that
$$
\begin{aligned}
|\mathcal{R}_S(\Lambda)-\mathcal{R}_T(\Lambda)|
&=\left|\int_0^\Lambda\bigl(\mathcal{N}_S(t)-\mathcal{N}_T(t)\bigr)\,dt\right|\\
&\leq\int_0^\Lambda|\mathcal{N}_S(t)-\mathcal{N}_T(t)|\,dt\\
&\leq\int_0^\Lambda r\,dt\\
&=\Lambda r.
\end{aligned}
$$
\begin{flushright}
  $\Box$
\end{flushright}

This completes the proof of Proposition~\ref{prop:weyl-rank}. We will use this proposition to compare the eigenvalue counting functions and Riesz means of the original operator and its block approximation.

\subsection{Compact embedding}

$\mathbf{The\ proof\ of\ Proposition\ 1.1}$ If \textup{($\mathbf{W_0}$)} holds, for every $u\in C_c(V)$, it follows from the Hardy inequality \eqref{Hardy-inequality} that
\begin{equation}\label{bedding}
\begin{aligned}
\|u\|_{\ell^2_W}^2
&=
\sum_{x\in V}\mu(x)W(x)|u(x)|^2\leq
M\sum_{x\in V}\mu(x)H_o(x)|u(x)|^2 \leq CM\|u\|_{1}^2.
\end{aligned}
\end{equation}
By density, the continuous embedding extends from $C_c(V)$ to $\cH^1=\overline{C_c(V)}^{\,\|\cdot\|_{1}}$, so the embedding is continuous.

Besides, if \textup{($\mathbf{W_1}$)} holds, $\forall \varepsilon>0$, we can choose $n\in\mathbb{N}$ sufficiently large that $W(x)\leq\varepsilon H_o(x)$ for every $x\notin B_n(0)$, where
$$
B_n(0)=\{x\in V:d_G(0,x)\leq n\}.
$$
Since $B_n(0)$ is finite and $H_o$ is nonnegative, there exists $C_n>0$ such that
\begin{equation}\label{C1}
W(x)\leq\varepsilon H_o(x),\quad x\in V\setminus B_n(0),
\qquad
W(x)\leq C_n H_o(x),\quad x\in B_n(0).
\end{equation}
Choose $M=\max\{\varepsilon, C_n\}$, we have
$$
W(x)<MH_o(x),\ \forall x\in V,
$$
so the embedding is also continuous.

Next, let $\{u_n\}_{n\in \mathbb{N}}$ be an arbitrary bounded sequence in $\cH^1$, so there exists  a constant $C>0$ such that
\begin{equation}\label{C2}
\|u_n\|_{1}\leq C.
\end{equation}
By the continuous embedding of $\cH^1$ to $\ell^2_{W}$, there exists  a constant $C>0$ such that
$$
\|u_n\|_{\ell^2_{W}}^2\leq C.
$$
Thus $\{u_n(x)\}_{n\in \mathbb{N}}$ is bounded in $\ell^2_{W}$, there exists $u\in\ell^2_{W}$ such that
$$
u_n\rightharpoonup u\quad \text{in }\ell^2_{W},\qquad u_n(x)\rightarrow u(x)\quad\text{for every }x\in V \quad\text{as }n\rightarrow\infty.
$$
By Fatou's lemma, \eqref{C2} and the Hardy inequality \eqref{Hardy-inequality} we have
\begin{equation}\label{C6}
\sum_{x\in V}\mu(x)H_o(x)|u(x)|^2\leq\liminf_{n\to\infty}\sum_{x\in V}\mu(x)H_o(x)|u_n(x)|^2\leq C.
\end{equation}
It follows from \eqref{C1}, \eqref{C6} and the inequality $|a-b|^2\leq2|a|^2+2|b|^2$ we have
\begin{equation}\label{C4}
\begin{aligned}
\sum_{x\notin B_n(0)}\mu(x)W(x)|u_n(x)-u(x)|^2
&\leq
\varepsilon\sum_{x\in V}\mu(x)H_o(x)|u_n(x)-u(x)|^2\\
&<C\varepsilon.
\end{aligned}
\end{equation}
Besides, since $B_n(0)$ is finite, so we have
\begin{equation}\label{C5}
\sum_{x\in B_n(0)}\mu(x)W(x)|u_n(x)-u(x)|^2\rightarrow0
\end{equation}
It follows from \eqref{C4} and \eqref{C5} that  $u_n\rightarrow u$ in $\ell^2_{W}$ as $n\rightarrow\infty$, which leads to the compact embedding.\hfill$\Box$

\section{Eigenvalue sequences}

For proving the Theorem 1.2, we need the Riesz-Fr\'echet representation theorem that needs to be used, which is the theorem 5.5 in \cite{Brezis}.

\begin{lemma}\label{num:riesz}
Let $\cH^1$ be a Hilbert space whose inner product is linear in the first argument. For every bounded linear functional $F$ on $\cH^1$, there exists a unique $u\in\cH^1$ such that
$$
F(v)=(v,u)_{1}
\quad\text{for all }v\in\cH^1.
$$
Moreover, $\|F\|_{(\cH^1)^*}=\|u\|_{1}$.
\end{lemma}

Now we are in a position to show the basic properties of eigenvalues.\medskip

\noindent{\bf The  proof  of Theorem 1.2 $(i)$. } Set
\begin{equation}\label{x1}
\lambda_1=\inf_{\substack{u\in \cH^1\\ \|u\|_{\ell_W^2}=1}}\|u\|_{1}^2.
\end{equation}
Since $\operatorname{supp}W\neq\varnothing$, we can choose $x_0\in \operatorname{supp}W$ and take
$$
u=\frac{\mathbf{1}_{\{x_0\}}}{\sqrt{\mu(x_0)W(x_0)}}\in C_c(V),
\quad \|u\|_{\ell_W^2}=1.
$$
By Proposition 1.1, the embedding of $\cH^1$ to $\ell_W^2$ is continuous \eqref{bedding}, so there exists $C>0$ such that
$$
1=\|u\|_{\ell_W^2}^2\le C\|u\|_{1}^2,
$$
it follows from \eqref{x1} that $\lambda_1\ge C^{-1}>0$.

Let $\{u_n\}\subset\cH^1$ be a minimizing sequence of $\lambda_1$ satisfying
$$
\|u_n\|_{\ell_W^2}=1,
\qquad
\|u_n\|_{1}^2\rightarrow\lambda_1\quad\text{as }n\rightarrow\infty.
$$
Since $\{u_n\}$ is bounded in $\cH^1$, there is a weak convergence subsequence still represented by $u_n$, which satisfies
$$
u_n\rightharpoonup \xi_1\quad\text{as }n\rightarrow\infty,
$$
for $\xi_1\in \cH^1$. It follows from the compact embedding of $\cH^1$ to $\ell_W^2$ that $\|\xi_1\|_{\ell_W^2}=1$. By weak lower semicontinuity,
$$
\lambda_1\le \|\xi_1\|_{1}^2\le \liminf_{n\to\infty}\|u_n\|_{1}^2=\lambda_1.
$$
Thus $\|\xi_1\|_{1}^2=\lambda_1$.

For $k\ge2$, suppose that $\xi_1,\ldots,\xi_{k-1}$ have been constructed and satisfy $(\xi_i,\xi_j)_{\ell_W^2}=\delta_{ij}$.
Set
$$
\lambda_k=\inf_{\substack{u\in \cH^1_k\\ \|u\|_{\ell_W^2}=1}}\|u\|_{1}^2,
$$
where we recall
$$
\cH^1_k=\{u\in\cH^1:(u,v)_{\ell_W^2}=0\text{ for every }v\in \operatorname{span}\{\xi_1,\ldots,\xi_{k-1}\}\}.
$$
If $k\le\dim\ell_W^2$, the density of the restrictions of $C_c(V)$ in $\ell_W^2$ allows us to choose $v\in C_c(V)$ with
$$
v\notin\operatorname{span}\{\xi_1,\ldots,\xi_{k-1}\}.
$$
Let
$$
v_0=v-\sum_{i=1}^{k-1}(v,\xi_i)_{\ell_W^2} \xi_i.
$$
Then $\|v_0\|_{\ell_W^2}>0$, and for every $j<k$,
$$
(v_0,\xi_j)_{\ell_W^2}=(v,\xi_j)_{\ell_W^2}-\sum_{i=1}^{k-1}(v,\xi_i)_{\ell_W^2}\delta_{ij}=0.
$$
Consequently, $\|\frac{v_0}{\|v_0\|_{\ell_W^2}}\|_{\ell_W^2}=1$. Choose a minimizing sequence $\{u_n\}\subset \cH^1_k$ with
$$
\|u_n\|_{\ell_W^2}=1,
\qquad
\|u_n\|_{1}^2\rightarrow\lambda_k\quad\text{as }n\rightarrow\infty.
$$
Along a subsequence $\{u_n\}\subset\cH^1$,
$$
u_n\rightharpoonup \xi_k,
\qquad
\|u_n-\xi_k\|_{\ell_W^2}\rightarrow0\quad \text{ as }n\rightarrow\infty.
$$
For each $i<k$, since $(u_n,\xi_i)_{\ell_W^2}=0$, we have
$$
0\le|(\xi_k,\xi_i)_{\ell_W^2}|=|(\xi_k-u_n,\xi_i)_{\ell_W^2}|\le \|\xi_k-u_n\|_{\ell_W^2}\|\xi_i\|_{\ell_W^2}\rightarrow0\quad \text{ as }n\rightarrow\infty.
$$
Thus $\xi_k\in \cH^1_k$. Moreover,
$$
0\le\bigl|\|\xi_k\|_{\ell_W^2}-1\bigr|=\bigl|\|\xi_k\|_{\ell_W^2}-\|u_n\|_{\ell_W^2}\bigr|\le \|\xi_k-u_n\|_{\ell_W^2}\rightarrow0\quad \text{ as }n\rightarrow\infty.
$$
Hence $\|\xi_k\|_{\ell_W^2}=1$, and
$$
\lambda_k\le \|\xi_k\|_{1}^2\le \liminf_{n\to\infty}\|u_n\|_{1}^2=\lambda_k.
$$

Next we verify the eigenvalue equation at each step. For $v\in \cH^1_k$,
$$
0=\left.\frac{d}{dt}\right|_{t=0}\frac{\|\xi_k+tv\|_{1}^2}{\|\xi_k+tv\|_{\ell_W^2}^2}=2\bigl((\xi_k,v)_{1}-\lambda_k(\xi_k,v)_{\ell_W^2}\bigr).
$$
we obtain
$$
(\xi_k,v)_{1}=\lambda_k(\xi_k,v)_{\ell_W^2},\qquad v\in \cH^1_k.
$$
For $k=1$, this is the eigenvalue equation on $\cH^1$. For $k\ge2$, the equations for the previously constructed eigenfunctions give
$$
(\xi_k,\xi_i)_{1}=\lambda_i(\xi_k,\xi_i)_{\ell_W^2}=0,
\quad i<k.
$$
For any $v\in \cH^1$,
$$
v-\sum_{i=1}^{k-1}(v,\xi_i)_{\ell_W^2} \xi_i\in \cH^1_k.
$$
The equation for $\xi_k$ therefore holds for every $v\in \cH^1$. For each $x\in V$, testing the weak eigenvalue equation with $\mathbf{1}_{\{x\}}\in C_c(V)\subset\cH^1$ yields $-\Delta \xi_k(x)=\lambda_k W(x)\xi_k(x)$. Finally, $\cH^1_{k+1}\subset \cH^1_k$ implies $\lambda_k\le\lambda_{k+1}$ whenever both are defined.

(ii) If $\#\operatorname{supp}W=m<\infty$, then $\dim\ell_W^2=m$. The construction yields $\xi_1,\ldots,\xi_m$. Their restrictions form an orthonormal basis of $\ell_W^2$, so if $u\in \cH^1_{m+1}$ then $\|u\|_{\ell_W^2}=0$. Thus the construction stops at $m$.

If $\operatorname{supp}W$ is infinite, then $\dim\ell_W^2=\infty$, and the construction continues for every $k$. If $\{\lambda_k\}$ were bounded, then
$$
\sup_k\|\xi_k\|_{1}^2=\sup_k\lambda_k<\infty.
$$
Compactness would give a subsequence converging in $\ell_W^2$.
However,
$$
\|\xi_i-\xi_j\|_{\ell_W^2}^2=2,\qquad i\ne j.
$$
Thus $\{\lambda_k\}$ is unbounded, and its monotonicity gives $\lambda_k\to+\infty$ as $k\rightarrow\infty$.

\medskip
\noindent\textup{(iii)} By construction and the eigenvalue equations,
$$
(\xi_i,\xi_j)_{\ell_W^2}=\delta_{ij},
\qquad
(\xi_i,\xi_j)_{1}=\lambda_i\delta_{ij}.
$$
Completeness in $\ell_W^2$ follows from dimension when $\operatorname{supp}W$ is finite. Suppose now that $\operatorname{supp}W$ is infinite and that $f\in\ell_W^2$ satisfies $(f,\xi_k)_{\ell_W^2}=0$ for every $k$. By the Riesz representation theorem i.e. Lemma~\ref{num:riesz}, there is $u\in \cH^1$ such that
$$
(u,v)_{1}=(f,v)_{\ell_W^2},\quad v\in \cH^1.
$$
For every $k$,
$$
\lambda_k(u,\xi_k)_{\ell_W^2}=(u,\xi_k)_{1}=(f,\xi_k)_{\ell_W^2}=0,
$$
so $u\in \cH^1_k$ for every $k$. If $\|u\|_{\ell_W^2}>0$, then
$$
\lambda_k\le \frac{\|u\|_{1}^2}{\|u\|_{\ell_W^2}^2}
\qquad\text{for every }k,
$$
which contradicts $\lambda_k\to+\infty$ as $k\rightarrow\infty$. Hence
$$
\|u\|_{\ell_W^2}=0,
\qquad
\|u\|_{1}^2=(f,u)_{\ell_W^2}=0,
$$
which implies $(f,v)_{\ell_W^2}=0$ for every $v\in \cH^1$ and $f=0$. Thus, the restrictions of $\{\xi_k\}$ form an orthonormal basis of $\ell_W^2$.
\begin{flushright}
$\Box$
\end{flushright}
\begin{corollary}
Assume that $\operatorname{supp}W$ is finite. Suppose that $-\Delta$ admits a Green function $\mathcal{G}$ for the energy space $\cH^1$, normalized so that, for every $y\in V$, $\mathcal{G}(\cdot,y)\in\cH^1$ and
$$
(\mathcal{G}(\cdot,y),v)_{1}=\mu(y)v(y),
\qquad v\in\cH^1.
$$
Then each eigenfunction $\xi_k$ from Theorem 1.2 satisfies
$$
\xi_k(x)=\lambda_k\sum_{y\in \operatorname{supp}W}\mathcal{G}(x,y)W(y)\xi_k(y),
\qquad x\in V.
$$
In particular, if $\mathcal{G}(x,y)\to 0$ as $x\to\infty$ for each $y\in \operatorname{supp}W$, then $\xi_k(x)\to 0$ as $x\to\infty$.
\end{corollary}

\begin{proof}
Set $\overline{u}_k=\lambda_k\sum_{y\in \operatorname{supp}W}W(y)\xi_k(y)\mathcal{G}(\cdot,y)$. Since $\operatorname{supp}W$ is finite, $\overline{u}_k\in\cH^1$. For every $v\in\cH^1$, the defining identity for $\mathcal{G}$ gives
$$
(\overline{u}_k,v)_{1}=\lambda_k\sum_{y\in \operatorname{supp}W}\mu(y)W(y)\xi_k(y)v(y)=\lambda_k(\xi_k,v)_{\ell_W^2}.
$$
Taking $v=\xi_k-\overline{u}_k$ gives $\|\xi_k-\overline{u}_k\|_1^2=0$. Hence $\xi_k=\overline{u}_k$.
\end{proof}

\section{Weyl Limit}
We first prove the asymptotic formula for the model weight
$$
W_0(x)=(1+|x|)^{-\alpha}[\log(e+|x|)]^\beta,
\qquad x\in\mathbb{Z}^d.
$$
In particular, $W_0(0)=1$. Until the comparison with the general
weight in the proof of Theorem~1.3, we write $W$ for $W_0$. We put
$$
w(r)=(1+r)^{-\alpha}[\log(e+r)]^\beta,\qquad r\geq0,
$$
where $\alpha>2$ and $\beta\in\mathbb{R}$. Since $ W(x) (1+|x|)^{2}\to0$ as $|x|\to\infty$, the Hardy inequality \eqref{Htineq} and condition (\ref{cond tilde W}) give the compact embedding required for Theorem 1.2.

\subsection{Finite volume approximation}

 Here we use the normalization $Au=\lambda Wu$. If the weighted square norm of Theorem 1.2 is used with $\mu(x)=2d$ and the same energy, the eigenvalues are divided by $2d$. Let $\{y_j\}_{j\geq1}$ be an enumeration of $\mathbb{Z}^d$ without repetition. For $L\in\mathbb{N}$ and $j\geq1$, define
$$
Q_{j,L}=x_j+[0,L-1]_\mathbb{Z}^d,\quad x_j=Ly_j
$$
where $[0,L-1]_{\mathbb{Z}}=[0,L-1]\cap\mathbb{Z}=\{0,1,\ldots,L-1\}$ and $[0,L-1]_{\mathbb{Z}}^d$ denotes $$\underbrace{[0,L-1]_{\mathbb{Z}}\times\cdots\times[0,L-1]_{\mathbb{Z}}}_d.$$
Let $\delta Q_{j,L}$ denote the inner vertex boundary of $Q_{j,L}$, having a nearest neighbor outside $Q_{j,L}$. Then
\begin{equation}\label{bounded}
|\delta Q_{j,L}|
\leq 2dL^{d-1}.
\end{equation}

For $\Lambda>1$, let
$$
a_\Lambda=\alpha^{-\frac{\beta}{\alpha}}\Lambda^{\frac{1}{\alpha}}(\log\Lambda)^{\frac{\beta}{\alpha}},
\qquad S_\Lambda=\Lambda a_\Lambda^d,
$$
where $\alpha>2$ and $\beta\in \mathbb{R}$. Then $a_\Lambda\to\infty$ as $\Lambda\rightarrow\infty$ and
$$
\Lambda a_\Lambda^{-\alpha}(\log a_\Lambda)^\beta\rightarrow1\quad\text{as }\Lambda\rightarrow\infty,
$$
because $\log a_\Lambda\sim\alpha^{-1}\log\Lambda$ as $\Lambda\rightarrow\infty$. Fix $M>0$ and set
$$
R=Ma_\Lambda,
\qquad
B_R=\{x\in\mathbb{Z}^d:|x|<R\}.
$$
Choose an integer $L=L_\Lambda$ such that
$$
L_\Lambda\rightarrow\infty,
\qquad
\frac{L_\Lambda}{a_\Lambda}\rightarrow0\quad\rm{as}\ \Lambda\rightarrow\infty.
$$
Define
$$
J_{R,L}=\{j\geq1:Q_{j,L}\subseteq B_R\}
$$
and
$$
D_{R,L}=B_R\setminus\bigcup_{j\in J_{R,L}}Q_{j,L}.
$$
Thus $J_{R,L}$ indexes the cubes entirely contained in $B_R$, while $D_{R,L}$ consists of the remaining vertices.

For $u$ on $B_R$, extend $u$ by zero to $\mathbb Z^d\setminus B_R$. The Dirichlet restriction $A_{B_R}$ is obtained by applying the lattice Laplacian to this extension and restricting the result to $B_R$. Explicitly, the Dirichlet eigenvalue problem on $B_R$ is
$$
\begin{cases}
A_{B_R}u(x)=\lambda W(x)u(x),
& x\in B_R,\\
u(x)=0, & x\in\partial B_R,
\end{cases}
$$
and
$$
A_{B_R}u(x)
=
2d\,u(x)-\sum_{\substack{y\in B_R\\ |y-x|=1}}u(y).
$$

For $u$ on $Q_{j,L}$, let $\widetilde u$ be its periodic
extension to $\mathbb Z^d$, determined by $\widetilde u=u$ on
$Q_{j,L}$ and
$$
\widetilde u(x+Le_i)=\widetilde u(x),
\qquad x\in\mathbb Z^d,\quad i=1,\ldots,d,
$$
where $e_i$ is the $i$-th coordinate unit vector. The periodic eigenvalue problem on $Q_{j,L}$ is
$$
\begin{cases}
A_{Q_{j,L}}^{\mathrm{per}}u(x)=\lambda W(x)u(x),
& x\in Q_{j,L},\\
\widetilde u(x+Le_i)=\widetilde u(x),
& x\in\mathbb{Z}^d,\quad i=1,\ldots,d,\\
\widetilde u(x)=u(x), & x\in Q_{j,L},
\end{cases}
$$
where
$$
A_{Q_{j,L}}^{\mathrm{per}}u(x)
=
\sum_{i=1}^d
\bigl(2u(x)-\widetilde u(x+e_i)-\widetilde u(x-e_i)\bigr),
\qquad x\in Q_{j,L}.
$$

For example, if $x$ lies on the upper face in direction $i$, then $\widetilde u(x+e_i)=u(x-(L-1)e_i)$, the value on the opposite face. The lower face is treated by the same periodic extension.

We firstly compare the finite volume Laplacian with the direct sum of the periodic block Laplacians:

\begin{lemma}\label{num:block-rank}
For each fixed $M>0$,
$$
\left|\mathcal{R}_{A_{B_R}}(\Lambda)-\sum_{j\in J_{R,L}}\mathcal{R}_{A_{Q_{j,L}}^{\mathrm{per}}}(\Lambda)\right|=o\!\left(S_\Lambda\right)\quad \rm{as}\ \Lambda\rightarrow\infty.
$$
\end{lemma}

\begin{proof}
We first estimate the number of remaining vertices. For any $x\in D_{R,L}$, the unique partition cube containing $x$ is not entirely contained in $B_R$. Hence there exists a vertex $y$ in the same cube such that $y\notin B_R$. Since $x$ and $y$ are in the same cube, so we have
$$
|x-y|\leq\sqrt{d}\,L,
$$
we obtain
$$
|x|
\geq |y|-|x-y|
\geq R-\sqrt{d}\,L.
$$
Therefore,
$$
D_{R,L}
\subseteq
\left\{
x\in\mathbb{Z}^d:
R-\sqrt{d}\,L\leq |x|<R
\right\}.
$$
Because $R=Ma_{\Lambda}$ and $L=o(R)$ for fixed $M$, a lattice shell
estimate gives
\begin{equation}\label{D1}
|D_{R,L}|
\leq C_dR^{d-1}L.
\end{equation}
Next, let
$$
n=|J_{R,L}|.
$$
The complete cubes are disjoint, and each contains
$L^d$ vertices. Consequently,
$$
nL^d\leq |B_R|\leq C_dR^d\quad \text{and}\quad n\leq C_d\left(\frac{R}{L}\right)^d.
$$
It follows from $|\delta Q_{j,L}|\leq 2dL^{d-1}$ in \eqref{bounded} that
\begin{equation}\label{D2}
\sum_{j\in J_{R,L}}|\delta Q_{j,L}|
\leq
2dnL^{d-1}
\leq
C_d\frac{R^d}{L}.
\end{equation}

To compare matrices of the same size, use the decomposition
$$
\ell^2(B_R)=\left(\bigoplus_{j\in J_{R,L}}\ell^2(Q_{j,L})\right)\oplus\ell^2(D_{R,L})
$$
and define
$$
T=A_{B_R},
\qquad
\widetilde T=\left(\bigoplus_{j\in J_{R,L}}A_{Q_{j,L}}^{\mathrm{per}}\right)\oplus 0_{D_{R,L}},
$$
where $0_{D_{R,L}}$ is the zero matrix on
$\ell^2(D_{R,L})$.

Set
$$
E=D_{R,L}\cup\bigcup_{j\in J_{R,L}}\delta Q_{j,L}.
$$
At every vertex outside $E$, the corresponding rows of $T$ and $\widetilde T$ agree: the vertex lies in the interior of a complete cube, where the original and periodic Laplacians have identical matrix entries. Thus by the \eqref{D1} and \eqref{D2} we have
$$
\operatorname{rank}(T-\widetilde T)
\leq |E|
\leq
C_d\left(
\frac{R^d}{L}+R^{d-1}L
\right).
$$

Both matrices are nonnegative. The Riesz mean of $\widetilde T$ equals $\sum_{j\in J_{R,L}}\mathcal{R}_{A_{Q_{j,L}}^{\mathrm{per}}}(\Lambda)+\Lambda|D_{R,L}|$, since the zero block on $D_{R,L}$ contributes $\Lambda$ per vertex. By Proposition~\ref{prop:weyl-rank}(iii), \eqref{D1}, $R=Ma_\Lambda$ and the definitions of $a_\Lambda$ and $S_\Lambda$ we have
$$
\frac{1}{S_\Lambda}\left|\mathcal{R}_{A_{B_R}}(\Lambda)-\sum_{j\in J_{R,L}}\mathcal{R}_{A_{Q_{j,L}}^{\mathrm{per}}}(\Lambda)\right|
\leq
C_d\left(\frac{M^d}{L}+M^{d-1}\frac{L}{a_\Lambda}\right)\rightarrow0\ \rm{as} \ \Lambda\rightarrow\infty,
$$
for fixed $M$, where the term $\Lambda|D_{R,L}|$ is absorbed by \eqref{D1}.
\end{proof}

Lemma~\ref{num:block-rank} shows that the Riesz mean of the Dirichlet Laplacian on $B_R$ can be approximated by the sum of the Riesz means of the periodic Laplacians on the cubes fully contained in $B_R$, with an error of $o(S_\Lambda)$ as $\Lambda\to\infty$ for each fixed $M>0$. We next show that replacing the weight on each cube fully contained in $B_R$ by its value at a representative vertex gives a total error of $o(S_\Lambda)$ in the sum of the block Riesz means as $\Lambda\to\infty$, for each fixed $M>0$.

For each cube
$$
Q_{j,L}=x_j+[0,L-1]_\mathbb{Z}^d,
$$
define
$$
r_j=|x_j|,
\qquad
W_j=W(x_j),
$$
and
$$
\delta_{j,L}
=
\sup_{x\in Q_{j,L}}
\left|\frac{W(x)}{W_j}-1\right|.
$$

Write
$$
W_{Q_{j,L}}=\operatorname{diag}\bigl(W(x)\bigr)_{x\in Q_{j,L}}.
$$
Define the variable weight and frozen weight matrices by
$$
T_j^W=W_{Q_{j,L}}^{-\frac{1}{2}}A_{Q_{j,L}}^{\mathrm{per}}W_{Q_{j,L}}^{-\frac{1}{2}},
\qquad
T_j^{W_j}=W_j^{-1}A_{Q_{j,L}}^{\mathrm{per}}.
$$
Denote their eigenvalues in nondecreasing order by
$$
\lambda_{n,j}^{W},
\qquad
\lambda_{n,j}^{W_j},
\qquad n=1,\ldots,L^d.
$$
Their Riesz means are
$$
\mathcal{R}_j^W(\Lambda)=\sum_{n=1}^{L^d}(\Lambda-\lambda_{n,j}^{W})_+,
\qquad
\mathcal{R}_j^{W_j}(\Lambda)=\sum_{n=1}^{L^d}(\Lambda-\lambda_{n,j}^{W_j})_+.
$$

\begin{lemma}\label{num:frozen}
If $r_j\geq2\sqrt{d}\,L$, then
$$
\delta_{j,L}
\leq C_{d,\alpha,\beta}\frac{L}{1+r_j}.
$$
When $\delta_{j,L}\leq\frac{1}{2}$, we have
$$
\frac{\lambda_{n,j}^{W_j}}{1+\delta_{j,L}}
\leq
\lambda_{n,j}^{W}
\leq
\frac{\lambda_{n,j}^{W_j}}{1-\delta_{j,L}},
\qquad n=1,\ldots,L^d.
$$
Moreover, for every $\Lambda\geq0$,
$$
E_j(\Lambda):=\left|\mathcal{R}_j^W(\Lambda)-\mathcal{R}_j^{W_j}(\Lambda)\right|\leq2\delta_{j,L}\Lambda L^d,
$$
and
$$
\sum_{j\in J_{R,L}}E_j(\Lambda)=o\!\left(S_\Lambda\right)\quad \rm{as}\ \Lambda\rightarrow\infty.
$$
\end{lemma}

\begin{proof}
Step 1: For $x\in Q_{j,L}$, write
$$
x=x_j+(n_1,\ldots,n_d),
\qquad 0\leq n_i\leq L-1.
$$
Then
\begin{equation}\label{r4}
||x|-|x_j||\le|x-x_j|
=
\left(\sum_{i=1}^d n_i^2\right)^{\frac{1}{2}}
\leq\sqrt{d}(L-1)
\leq\sqrt{d}\,L.
\end{equation}
It follows from the assumption $r_j\geq2\sqrt{d}\,L$ that we have
$$
|x|
\geq |x_j|-|x-x_j|
\geq r_j-\sqrt{d}\,L
\geq\frac{r_j}{2}
$$
and
\begin{equation}\label{rrr1}
\frac{L}{1+r_j}
\leq \frac{L}{r_j}
\leq \frac{1}{2\sqrt{d}},
\qquad
\frac{1}{1+r}\leq \frac{2}{1+r_j}.
\end{equation}

For the radial profile $w(r)=(1+r)^{-\alpha}[\log(e+r)]^\beta$, we have
$$
\log w(r)=-\alpha\log(1+r)+\beta\log\log(e+r).
$$
Differentiating with respect to $r$ gives
$$
\frac{d}{dr}\log w(r)
=
-\frac{\alpha}{1+r}
+
\frac{\beta}{(e+r)\log(e+r)}.
$$
By the triangle inequality,
$$
\left|\frac{d}{dr}\log w(r)\right|
\leq
\frac{\alpha}{1+r}
+
\frac{|\beta|}{(e+r)\log(e+r)}.
$$
For every $r\geq 0$, we have $e+r\geq 1+r$ and $\log(e+r)\geq 1$. Hence
$$
(e+r)\log(e+r)\geq 1+r.
$$
Then we have
\begin{equation}\label{w1}
\left|\frac{d}{dr}\log w(r)\right|\leq \frac{\alpha+|\beta|}{1+r}.
\end{equation}
By the fundamental theorem of calculus, \eqref{r4}, \eqref{rrr1} and \eqref{w1} we have
\begin{equation}\label{diflog}
\begin{aligned}
\left|\log\frac{W(x)}{W_j}\right|
&=
\left|\int_{r_j}^{|x|}(\log w)'(r)\,dr\right|\\
&\leq
\frac{2(\alpha+|\beta|)}{1+r_j}\,\bigl||x|-r_j\bigr|\\
&\leq
2\sqrt{d}\,(\alpha+|\beta|)\frac{L}{1+r_j}\\
&\le\alpha+|\beta|.
\end{aligned}
\end{equation}

Set $s=\log\frac{W(x)}{W_j}$, then we have $|s|\leq\alpha+|\beta|$ uniformly in $j$, $L$, and $x$. Using $|e^s-1|\leq e^{|s|}|s|$ and \eqref{diflog}, we obtain
$$
\begin{aligned}
\left|\frac{W(x)}{W_j}-1\right|
&=|e^s-1|\\
&\leq e^{|s|}|s|\\
&\leq 2\sqrt{d}\,(\alpha+|\beta|)e^{\alpha+|\beta|}\frac{L}{1+r_j}.
\end{aligned}
$$
Taking the supremum over $x\in Q_{j,L}$ yields
$$
\delta_{j,L}
=
\sup_{x\in Q_{j,L}}\left|\frac{W(x)}{W_j}-1\right|
\leq C_{d,\alpha,\beta}\frac{L}{1+r_j},
$$
where $C_{d,\alpha,\beta}=2\sqrt{d}\,(\alpha+|\beta|)e^{\alpha+|\beta|}$
\\
\\
Step 2: By the definition of $\delta_{j,L}$, we have
\begin{equation}\label{W3}
(1-\delta_{j,L})W_j
\leq W(x)
\leq(1+\delta_{j,L})W_j,
\qquad x\in Q_{j,L}.
\end{equation}
For $u\neq0$, define the generalized Rayleigh quotients in $Q_{j,L}$
\begin{equation}\label{Q1}
Q_j^W(u)=\frac{u^*A_{Q_{j,L}}^{\mathrm{per}}u}{\displaystyle\sum_{x\in Q_{j,L}}W(x)|u(x)|^2},
\end{equation}
and
\begin{equation}\label{Q2}
Q_j^{W_j}(u)=\frac{u^*A_{Q_{j,L}}^{\mathrm{per}}u}{\displaystyle W_j\sum_{x\in Q_{j,L}}|u(x)|^2}.
\end{equation}
By the \eqref{W3} we have
\begin{equation}\label{sumw}
(1-\delta_{j,L})W_j\sum_{x\in Q_{j,L}}|u(x)|^2\leq\sum_{x\in Q_{j,L}}W(x)|u(x)|^2\leq(1+\delta_{j,L})W_j\sum_{x\in Q_{j,L}}|u(x)|^2.
\end{equation}
Since $A_{Q_{j,L}}^{\mathrm{per}}$ is nonnegative, if $\delta_{j,L}\leq\frac{1}{2}$ it follows from \eqref{Q1}, \eqref{Q2} and \eqref{sumw} that
$$
\frac{Q_j^{W_j}(u)}{1+\delta_{j,L}}
\leq Q_j^W(u)
\leq
\frac{Q_j^{W_j}(u)}{1-\delta_{j,L}}.
$$
By the Lemma~\ref{num:minmax} we have
$$
\frac{\lambda_{n,j}^{W_j}}{1+\delta_{j,L}}
\leq\lambda_{n,j}^{W}
\leq
\frac{\lambda_{n,j}^{W_j}}{1-\delta_{j,L}}.
$$
\\
Step 3: Fix $n$, for the sake of simplicity, let
$$
p=\lambda_{n,j}^{W_j},
\qquad
q=\lambda_{n,j}^{W},
\qquad
\delta=\delta_{j,L}.
$$
Then by step 2 we have
\begin{equation}\label{simplicity}
\frac{p}{1+\delta}
\leq q
\leq\frac{p}{1-\delta},
\qquad 0\leq\delta\leq\frac{1}{2}.
\end{equation}
We claim that
$$
|(\Lambda-p)_+-(\Lambda-q)_+|
\leq2\delta\Lambda.
$$

If $\min\{p,q\}\geq\Lambda$, the claim is clearly valid. Otherwise, consider the following two cases.

If $q\geq p$, then $p<\Lambda$, it follows from \eqref{simplicity} that
$$
0\leq q-p
\leq\frac{p}{1-\delta}-p
=
\frac{\delta}{1-\delta}p
\leq2\delta\Lambda,
$$
so we have
$$
|(\Lambda-p)_+-(\Lambda-q)_+|
\leq|p-q|
\leq2\delta\Lambda.
$$

If $p>q$, then $q<\Lambda$.
The inequality $\frac{p}{(1+\delta)}\leq q$ in \eqref{simplicity} gives
$$
p\leq(1+\delta)q,
$$
and hence
$$
0<p-q\leq\delta q\leq\delta\Lambda.
$$
so we have
$$
|(\Lambda-p)_+-(\Lambda-q)_+|
\leq|p-q|
\leq2\delta\Lambda.
$$

Summing over $n=1,\ldots,L^d$, we obtain
$$
\begin{aligned}
E_j(\Lambda)&=\left|\mathcal{R}_j^W(\Lambda)-\mathcal{R}_j^{W_j}(\Lambda)\right|\\
&\leq
\sum_{n=1}^{L^d}
\left|
(\Lambda-\lambda_{n,j}^{W})_+
-
(\Lambda-\lambda_{n,j}^{W_j})_+
\right|\\
&\leq
2\delta_{j,L}\Lambda L^d.
\end{aligned}
$$
\\
Step 4: Write $a=a_\Lambda$ and $L=L_{\Lambda}$. Fix $\varepsilon\in(0,1)$. We split the sum into the regions $r_j\leq\varepsilon a$ and $r_j>\varepsilon a$.

For the inner region, $r_j=L|y_j|$ and the enumeration has no repetitions, so
\begin{equation}\label{J1}
\#\{j\in J_{R,L}:r_j\leq\varepsilon a\}
\leq
C_d\left(1+\frac{\varepsilon a}{L}\right)^d.
\end{equation}
Both local matrices are nonnegative and have $L^d$
eigenvalues, so
$$
0\leq \mathcal{R}_j^W(\Lambda)\leq\Lambda L^d,
\qquad
0\leq \mathcal{R}_j^{W_j}(\Lambda)\leq\Lambda L^d.
$$
Consequently, without requiring any smallness assumption
on $\delta_{j,L}$,
$$
E_j(\Lambda)\leq\Lambda L^d.
$$
It follows from \eqref{J1} we have
\begin{equation}\label{in}
\sum_{\substack{j\in J_{R,L}\\r_j\leq\varepsilon a}}E_j(\Lambda)\leq C_d\Lambda L^d\left(1+\frac{\varepsilon a}{L}\right)^d=C_d\Lambda(L+\varepsilon a)^d.
\end{equation}

For the outer region, the condition $\frac{L}{a}\to0$ ensures that, for sufficiently large $\Lambda$,
$$
r_j>\varepsilon a\geq2\sqrt{d}\,L.
$$
Step~1 then gives the uniform estimate
\begin{equation}\label{small}
\delta_{j,L}
\leq C_{d,\alpha,\beta}\frac{L}{1+r_j}
\leq C_{d,\alpha,\beta}\frac{L}{\varepsilon a}
\rightarrow0\quad \rm{as}\ \Lambda\rightarrow\infty.
\end{equation}
In particular, $\delta_{j,L}\leq\frac{1}{2}$ throughout
this region for sufficiently large $\Lambda$.
By Step~3,
\begin{equation}\label{Ej1}
E_j(\Lambda)
\leq
C_{d,\alpha,\beta}\frac{L}{\varepsilon a}
\Lambda L^d.
\end{equation}
Since the complete cubes are disjoint and contained
in $B_R$,
\begin{equation}\label{Ej2}
\sum_{j\in J_{R,L}}L^d
\leq|B_R|
\leq C_dR^d
=
C_dM^da^d.
\end{equation}
It follows \eqref{Ej1} and \eqref{Ej2} that
\begin{equation}\label{out}
\sum_{\substack{j\in J_{R,L}\\r_j>\varepsilon a}}
E_j(\Lambda)
\leq
C_{d,\alpha,\beta}M^d
\frac{L}{\varepsilon a}S_\Lambda.
\end{equation}

Combining the \eqref{in} and \eqref{out}, we obtain
$$
0\leq
\frac{1}{S_\Lambda}\displaystyle\sum_{j\in J_{R,L}}E_j(\Lambda)
\leq
C_d\left(\frac{L}{a}+\varepsilon\right)^d
+
C_{d,\alpha,\beta}M^d\frac{L}{\varepsilon a}.
$$
First let $\Lambda\to\infty$ with $\varepsilon$
and $M$ fixed. Since $\frac{L}{a}\to0$ as $\Lambda\rightarrow\infty$,
$$
\limsup_{\Lambda\to\infty}
\frac{1}{S_\Lambda}\displaystyle\sum_{j\in J_{R,L}}E_j(\Lambda)
\leq C_d\varepsilon^d.
$$
Finally, letting $\varepsilon\rightarrow0$, we conclude
that
$$
\sum_{j\in J_{R,L}}E_j(\Lambda)
=
o(S_\Lambda)
=
o\!\left(S_\Lambda\right).
$$
\end{proof}

Combining the rank estimate from Lemma~\ref{num:block-rank} with Lemma~\ref{num:frozen}, we approximate the weighted Dirichlet Riesz mean on $B_R$ by the sum of the periodic
Riesz means with frozen weights over all cubes fully contained in $B_R$. The total error is $o(S_\Lambda)$ as $\Lambda\to\infty$ for each fixed $M>0$. Next, we diagonalize the periodic Laplacian in the Fourier basis and compare the resulting frequency sum with the integral defining $\rho$.

\begin{lemma}\label{num:periodic}
$$
\sum_{j\in J_{R,L}}\mathcal{R}_j^{W_j}(\Lambda)
=
\sum_{j\in J_{R,L}}L^dW_j^{-1}\rho(\Lambda W_j)
+
o\!\left(S_\Lambda\right)\quad \rm{as}\ \Lambda\rightarrow\infty.
$$
\end{lemma}

\begin{proof}
For each $m\in[0,L-1]_{\mathbb{Z}}^d$, define
$$
\varphi_m(x)=L^{-\frac{d}{2}}e^{2\pi i m\cdot\frac{x}{L}},
\qquad x\in Q_{j,L}.
$$
Let $e_r$ denote the $r$-th coordinate unit vector. Then we have
$$
\varphi_m(x+Le_r)=\varphi_m(x),\quad
\varphi_m(x+e_r)=e^{2\pi i \frac{m_r}{L}}\varphi_m(x),\quad
\varphi_m(x-e_r)=e^{-2\pi i \frac{m_r}{L}}\varphi_m(x).
$$
Substituting these identities into the periodic Laplacian gives
$$
A_{Q_{j,L}}^{\mathrm{per}}\varphi_m(x)
=
2\sum_{r=1}^{d}\left(1-\cos\frac{2\pi m_r}{L}\right)\varphi_m(x)
=
\sigma\left(\frac{2\pi m}{L}\right)\varphi_m(x).
$$
The functions $\varphi_m$ are an orthonormal basis of
$\ell^2(Q_{j,L})$ because
$$
\langle\varphi_m,\varphi_\ell\rangle=
\begin{cases}
1,&m=\ell,\\
0,&m\neq\ell.
\end{cases}
$$
Since there are $L^d$ such functions, they account for
all eigenvalues.

For the frozen weight eigenvalue problem
$$
A_{Q_{j,L}}^{\mathrm{per}}u=\lambda W_j u,
$$
the eigenvalues are therefore
$$
\lambda_{m,j}^{\mathrm{per},W_j}
=
W_j^{-1}\sigma\left(\frac{2\pi m}{L}\right).
$$
Consequently,
$$
\mathcal{R}_j^{W_j}(\Lambda)
=
\sum_{m}\left(\Lambda-W_j^{-1}
\sigma\left(\frac{2\pi m}{L}\right)\right)_+
=
W_j^{-1}\sum_{m}
\left(\Lambda W_j-\sigma\left(\frac{2\pi m}{L}\right)\right)_+.
$$
Define
$$
\rho_L(t)
=
\frac{1}{L^d}
\sum_{m}
\left(t-\sigma\left(\frac{2\pi m}{L}\right)\right)_+.
$$
Then we have
\begin{equation}\label{R1}
\mathcal{R}_j^{W_j}(\Lambda)
=
L^dW_j^{-1}\rho_L(\Lambda W_j).
\end{equation}

For $t\geq0$, let
$$
f_t(\theta)=(t-\sigma(\theta))_+.
$$
Thus
\begin{equation}\label{r3}
\rho(t)
=
\frac{1}{(2\pi)^d}
\int_{[0,2\pi)^d}f_t(\theta)\,d\theta,
\qquad
\rho_L(t)
=
\frac{1}{L^d}\sum_{m}f_t(\theta_m),
\end{equation}
where $\theta_m=\frac{2\pi m}{L}$. Let
\begin{equation}\label{Tm}
T_m
=
\prod_{r=1}^{d}
\left[
\frac{2\pi m_r}{L},
\frac{2\pi(m_r+1)}{L}
\right).
\end{equation}
Each $T_m$ has volume
\begin{equation}\label{T1}
|T_m|=\left(\frac{2\pi}{L}\right)^d.
\end{equation}
For every $\theta=(\theta_1,\cdots,\theta_d)\in T_m$,
$$
|\theta-\theta_m|^2
=
\sum_{r=1}^{d}
\left|\theta_r-\frac{2\pi m_r}{L}\right|^2
\leq
d\left(\frac{2\pi}{L}\right)^2,
$$
so
\begin{equation}\label{sitacha}
|\theta-\theta_m|
\leq\frac{2\pi\sqrt{d}}{L}.
\end{equation}
Moreover,
$$
\nabla\sigma(\theta)
=
(2\sin\theta_1,\ldots,2\sin\theta_d),
\qquad
|\nabla\sigma(\theta)|\leq2\sqrt{d},
$$
it follows from \eqref{sitacha} that
$$
|\sigma(\theta)-\sigma(\theta_m)|
\leq
2\sqrt{d}\,|\theta-\theta_m|
\leq\frac{4\pi d}{L}.
$$
It follows from $|(t-a)_+-(t-b)_+|\leq|a-b|$ and the definition of $f_t(\theta)=(t-\sigma(\theta))_+$, we obtain
\begin{equation}\label{f1}
|f_t(\theta)-f_t(\theta_m)|
\leq\frac{4\pi d}{L},
\end{equation}
uniformly in $t\geq0$. Using the \eqref{r3} and \eqref{Tm}, we may write
$$
\rho_L(t)
=
\frac{1}{(2\pi)^d}
\sum_{m}
\int_{T_m}f_t(\theta_m)\,d\theta.
$$
It follows \eqref{T1}, \eqref{f1} and $m\in[0,L-1]_{\mathbb{Z}}^d$ that
$$
\begin{aligned}
|\rho_L(t)-\rho(t)|
&\leq
\frac{1}{(2\pi)^d}
\sum_{m}
\int_{T_m}
|f_t(\theta_m)-f_t(\theta)|\,d\theta\\
&\leq
\frac{1}{(2\pi)^d}
\sum_{m}
|T_m|\frac{4\pi d}{L}\\
&=
\frac{4\pi d}{L}.
\end{aligned}
$$
Thus
\begin{equation}\label{r2}
\sup_{t\geq0}|\rho_L(t)-\rho(t)|
\leq\frac{4\pi d}{L}.
\end{equation}
By the \eqref{R1} and \eqref{r2} we have
\begin{equation}\label{R2}
\left|
\mathcal{R}_j^{W_j}(\Lambda)-L^dW_j^{-1}\rho(\Lambda W_j)\right|
=L^dW_j^{-1}|\rho_L(\Lambda W_j)-\rho(\Lambda W_j)|
\leq\frac{4\pi d}{L}L^dW_j^{-1}.
\end{equation}
Fix $M>0$ and set $R=Ma_\Lambda$, where
$$
a_\Lambda
=
\alpha^{-\frac{\beta}{\alpha}}
\Lambda^{\frac{1}{\alpha}}(\log\Lambda)^{\frac{\beta}{\alpha}}.
$$
We next show that there exists a constant $C_M>0$, independent of $\Lambda$, $L$ and $j$, such that
\begin{equation}\label{Wj2}
W_j^{-1}\leq C_M\Lambda
\end{equation}
for all representatives $x_j\in B_R$ and all sufficiently large
$\Lambda$.

Recall that $w(r)=(1+r)^{-\alpha}[\log(e+r)]^\beta$. Its
logarithmic derivative satisfies
$$
(1+r)\frac{d}{dr}\log w(r)
=
-\alpha
+
\beta\frac{1+r}{(e+r)\log(e+r)}
\rightarrow -\alpha<0\quad\text{as }r\to\infty.
$$
Since $w(r)>0$, there exists $r_0>0$ such that $w$ is decreasing on
$[r_0,\infty)$. Besides, by the definition of $a_\Lambda$ and $w(x)$
we have
$$
\Lambda w(Ma_\Lambda)
=
\Lambda a_\Lambda^{-\alpha}(\log a_\Lambda)^\beta
\left(M+\frac{1}{a_\Lambda}\right)^{-\alpha}
\left(\frac{\log(e+Ma_\Lambda)}{\log a_\Lambda}\right)^\beta
\rightarrow M^{-\alpha}\quad\text{as }\Lambda\rightarrow\infty.
$$
In particular, for all sufficiently large $\Lambda$, we have
$R>r_0$ and
$$
w(R)\geq \frac{M^{-\alpha}}{2\Lambda}.
$$

If $r_0\leq |x_j|<R$, the monotonicity of $w$ gives
$$
W_j=w(|x_j|)\geq w(R)\geq \frac{M^{-\alpha}}{2\Lambda},
$$
and hence $W_j^{-1}\leq 2M^\alpha\Lambda$.
For the remaining representatives, define
$$
C_0
=
\max_{\substack{x\in\mathbb Z^d\\ |x|<r_0}}W(x)^{-1}.
$$
This maximum is finite because the set is finite and $W(x)>0$ at every vertex, including the origin. Thus, if $|x_j|<r_0$ and $\Lambda\geq1$, then $W_j^{-1}\leq C_0\leq C_0\Lambda$. Taking $C_M=\max\{2M^\alpha,C_0\}$ proves the \eqref{Wj2}.

Besides, the cubes indexed by $J$ are disjoint and contained in $B_R$, so
\begin{equation}\label{Wj1}
\sum_{j\in J_{R,L}}L^d
=
\left|\bigcup_{j\in J_{R,L}}Q_{j,L}\right|
\leq|B_R|\leq C_dR^d.
\end{equation}

It follows from \eqref{R2}, \eqref{Wj2} and \eqref{Wj1} we have
$$
\left|\sum_{j\in J_{R,L}}\mathcal{R}_j^{W_j}(\Lambda)
-\sum_{j\in J_{R,L}}L^dW_j^{-1}\rho(\Lambda W_j)\right|
\leq C_{d,\alpha,\beta,M}\frac{S_\Lambda}{L}.
$$
Since $L\to\infty$ as $\Lambda\rightarrow\infty$, the last expression is $o(S_\Lambda)$. This proves the lemma.
\end{proof}

In the next lemma, we compare the frozen weight sum with the sum over vertices by estimating the cubes near the origin, the remaining complete cubes, and the boundary layer separately.

\begin{lemma}\label{num:lattice}
$$
\sum_{j\in J_{R,L}}L^dW_j^{-1}\rho(\Lambda W_j)
=
\sum_{x\in B_R}W(x)^{-1}\rho(\Lambda W(x))
+
o(S_\Lambda).
$$
\end{lemma}

\begin{proof}
By the lemma~\ref{num:symbol} we have $0\leq\rho(t)\leq t$, then
\begin{equation}\label{rr3}
0\leq W(x)^{-1}\rho(\Lambda W(x))\leq\Lambda,\quad 0\leq W_j^{-1}\rho(\Lambda W_j)\leq\Lambda.
\end{equation}
Each cube contains $L^d$ vertices, so
$$
\left|L^dW_j^{-1}\rho(\Lambda W_j)-\sum_{x\in Q_{j,L}}W(x)^{-1}\rho(\Lambda W(x))\right|\leq\sum_{x\in Q_{j,L}}|W_j^{-1}\rho(\Lambda W_j)-W(x)^{-1}\rho(\Lambda W(x))|\leq\Lambda L^d.
$$
Fix $\varepsilon\in(0,1)$. We next estimate separately the cubes with $|x_j|\leq\varepsilon a$ and those with $|x_j|>\varepsilon a$.

Since $|x_j|=L|y_j|$ and the enumeration has no repetitions, so we have
$$
\#\{j\in J_{R,L}:|x_j|\leq\varepsilon a\}
\leq
C_d\left(1+\frac{\varepsilon a}{L}\right)^d.
$$
Therefore,
\begin{equation}\label{10.1}
\begin{aligned}
&
\sum_{\substack{j\in J_{R,L}\\|x_j|\leq\varepsilon a}}
\left|
L^dW_j^{-1}\rho(\Lambda W_j)-\sum_{x\in Q_{j,L}}W(x)^{-1}\rho(\Lambda W(x))
\right|\\
&\qquad\leq
C_d\Lambda L^d
\left(1+\frac{\varepsilon a}{L}\right)^d\\
&\qquad=
C_d\Lambda(L+\varepsilon a)^d\\
&\qquad=
C_dS_\Lambda
\left(\frac{L}{a}+\varepsilon\right)^d.
\end{aligned}
\end{equation}

Secondly, for $|x_j|>\varepsilon a$, the condition $\frac{L}{a}\to0$ implies that, for sufficiently large $\Lambda$,
$$
|x_j|\geq2\sqrt{d}L,
$$
and
$$
\delta_{j,L}=\sup_{x\in Q_{j,L}}\left|\frac{W(x)}{W_j}-1\right|\leq\frac{1}{2}.
$$
Arguing as in Step 3 of the proof of Lemma~\ref{num:frozen} and using \eqref{small} we have
$$
|L^dW_j^{-1}\rho(\Lambda W_j)-\sum_{x\in Q_{j,L}}W(x)^{-1}\rho(\Lambda W(x))|\leq2\delta_{j,L}\Lambda L^d,
$$
and
\begin{equation}\label{10.2}
\begin{aligned}
&
\sum_{\substack{j\in J_{R,L}\\|x_j|>\varepsilon a}}
\left|
L^dW_j^{-1}\rho(\Lambda W_j)-\sum_{x\in Q_{j,L}}W(x)^{-1}\rho(\Lambda W(x))
\right|\\
&\qquad\leq
C_{d,\alpha,\beta}\frac{L}{\varepsilon a}
\Lambda\sum_{j\in J_{R,L}}L^d\\
&\qquad\leq
C_{d,\alpha,\beta}\frac{L}{\varepsilon a}\Lambda R^d\\
&\qquad=
C_{d,\alpha,\beta}M^d
\frac{L}{\varepsilon a}S_\Lambda,
\end{aligned}
\end{equation}
where we used $R=Ma_{\Lambda}$. Besides, by \eqref{D1} we have
$$
|D_{R,L}|\leq C_dR^{d-1}L,
$$
using $0\leq W(x)^{-1}\rho(\Lambda W(x))\leq\Lambda$ in \eqref{rr3}, we obtain
\begin{equation}\label{10.3}
0\leq\sum_{x\in D_{R,L}}W(x)^{-1}\rho(\Lambda W(x))\leq\Lambda|D_{R,L}|\leq C_d\Lambda R^{d-1}L=C_dM^{d-1}\frac{L}{a}S_\Lambda.
\end{equation}

Combining the \eqref{10.1}, \eqref{10.2} and \eqref{10.3} we have
$$
\begin{aligned}
~&\frac{1}{S_\Lambda}
\left|
\sum_{j\in J_{R,L}}L^dW_j^{-1}\rho(\Lambda W_j)
-
\sum_{x\in B_R}W(x)^{-1}\rho(\Lambda W(x))
\right|\\
&\qquad\leq
C_d\left(\frac{L}{a}+\varepsilon\right)^d
+
C_{d,\alpha,\beta}M^d\frac{L}{\varepsilon a}
+
C_dM^{d-1}\frac{L}{a}.
\end{aligned}
$$
First let $\Lambda\to\infty$, keeping $M$ and $\varepsilon$ fixed. Since $\frac{L}{a}\to0$ as $\Lambda\rightarrow\infty$, we have
$$
\limsup_{\Lambda\to\infty}
\frac{1}{S_\Lambda}
\left|
\sum_{j\in J_{R,L}}L^dW_j^{-1}\rho(\Lambda W_j)
-
\sum_{x\in B_R}W(x)^{-1}\rho(\Lambda W(x))
\right|
\leq C_d\varepsilon^d,
$$
then letting $\varepsilon\rightarrow0$ we have
$$
\sum_{j\in J_{R,L}}L^dW_j^{-1}\rho(\Lambda W_j)
=
\sum_{x\in B_R}W(x)^{-1}\rho(\Lambda W(x))
+
o(S_\Lambda).
$$
This completes the proof.
\end{proof}

\subsection{Asymptotics of the Riesz mean}

The main result of this subsection is the asymptotic formula for the full space Riesz mean $\mathcal{R}(\Lambda)=\sum_{k\geq1}(\Lambda-\lambda_k)_+$ of the eigenvalues $\{\lambda_k\}_{k\geq1}$ of \eqref{E2}.
\begin{weylproposition}\label{prop:riesz-asymptotic}
Define
$$
K_{d,\alpha}=\int_{\mathbb{R}^d}|y|^\alpha\rho(|y|^{-\alpha})\,dy,
$$
Then
$$
0<K_{d,\alpha}<\infty
$$
and
$$
\mathcal{R}(\Lambda)
\sim
K_{d,\alpha}S_\Lambda\quad\text{as }\Lambda\rightarrow\infty.
$$
\end{weylproposition}

The proof first reduces the full space Riesz mean to a lattice sum. We compare the Dirichlet Laplacian on $B_R$ with the direct sum of periodic block Laplacians by estimating the rank of their difference. On each complete block, we replace the variable weight $W(x)$ by the constant $W_j=W(x_j)$. Fourier diagonalization then gives a sum over the discrete frequencies of each block. Replacing these frequency sums by the integral defining $\rho$ gives
$$
\sum_{j\in J_{R,L}} L^d W_j^{-1}\rho(\Lambda W_j).
$$
We next replace this block sum by the vertex sum over $B_R$. Estimates for the exterior spectral contribution and the lattice sum tail give
$$
\mathcal{R}(\Lambda)
=
\sum_{x\in\mathbb{Z}^d}W(x)^{-1}\rho(\Lambda W(x))
+o(S_\Lambda)\quad\text{as }\Lambda\to\infty.
$$
Finally, we rescale the lattice sum by setting $y=\frac{x}{a_\Lambda}$. After division by $S_\Lambda$, the rescaled sum converges to the integral defining $K_{d,\alpha}$.

For $d\geq3$, we use the estimate that follows from Theorem 1.1 of \cite{BachPedraLakaev} by zero extension and the variational principle:
\begin{equation}\label{CLR}
\mathcal{N}_-(A_\Omega-V)
\leq C_d\sum_{x\in\Omega}V(x)^{\frac{d}{2}},
\qquad V\geq0,
\end{equation}
where $\Omega\subset\mathbb Z^d$ is finite, $A_\Omega$ is the Dirichlet restriction of the nonnegative lattice Laplacian, and $C_d$ is independent of $\Omega$ and $V$. Here $\mathcal{N}_-$ counts strictly negative eigenvalues, with multiplicity.

\begin{lemma}\label{num:full}
$$
\mathcal{R}(\Lambda)
=
\sum_{x\in\mathbb{Z}^d}
W(x)^{-1}\rho(\Lambda W(x))
+
o(S_\Lambda)\quad\text{as }\Lambda\rightarrow\infty.
$$
\end{lemma}

\begin{proof}
For $N>2R$, define
$$
\Omega_N=B_N,
\qquad
E_{R,N}=\Omega_N\setminus B_R.
$$
For any finite set $\Omega$, write
$$
T_\Omega
:=
W_\Omega^{-\frac{1}{2}}A_\Omega W_\Omega^{-\frac{1}{2}},
$$
and let $\mathcal{R}_\Omega$ denote its Riesz mean.

On
$$
\ell^2(\Omega_N)
=
\ell^2(B_R)\oplus\ell^2(E_{R,N}),
$$
define
$$
\widetilde T_N=T_{B_R}\oplus T_{E_{R,N}}.
$$

If $x\in B_R$, $y\in E_{R,N}$, and $|x-y|=1$,
then
$$
R-1\leq|x|<R,
\qquad
R\leq|y|<R+1.
$$
Thus all affected vertices lie in a shell containing
at most $C_dR^{d-1}$ vertices. Diagonal conjugation
by the positive weight does not change rank, so
$$
\operatorname{rank}(T_{\Omega_N}-\widetilde T_N)
\leq C_dR^{d-1}.
$$
By Proposition~\ref{prop:weyl-rank}(iii), we have
\begin{equation}\label{R-RBR-RER}
\mathcal{R}_{\Omega_N}(\Lambda)
\leq
\mathcal{R}_{B_R}(\Lambda)
+
\mathcal{R}_{E_{R,N}}(\Lambda)
+
C_d\Lambda R^{d-1}.
\end{equation}

For $t>0$,
$$
\#\{k:\lambda_k(T_{E_{R,N}})<t\}=\mathcal{N}_-(A_{E_{R,N}}-tW_{E_{R,N}}).
$$
Applying the \eqref{CLR} with $V=tW_{E_{R,N}}$ we have
\begin{equation}\label{NE}
\#\{k:\lambda_k(T_{E_{R,N}})<t\}\leq C_dt^{\frac{d}{2}}\sum_{x\in E_{R,N}}W(x)^{\frac{d}{2}}\leq C_dt^{\frac{d}{2}}\sum_{|x|\geq R}W(x)^{\frac{d}{2}}.
\end{equation}
So we have
\begin{equation}\label{RE}
\mathcal{R}_{E_{R,N}}(\Lambda)=\int_0^\Lambda \mathcal{N}_{T_{E_{R,N}}}(t)\,dt\leq C_d\Lambda^{1+\frac{d}{2}}\sum_{|x|\geq R}W(x)^{\frac{d}{2}}.
\end{equation}
This bound is independent of $N$. Let $N\rightarrow\infty$, by \eqref{R-RBR-RER} we have
\begin{equation}\label{R3}
0\leq \mathcal{R}(\Lambda)-\mathcal{R}_{B_R}(\Lambda)\leq C_d\Lambda R^{d-1}+C_d\Lambda^{1+\frac{d}{2}}\sum_{|x|\geq R}W(x)^{\frac{d}{2}}.
\end{equation}

Fix $\eta\in(0,\alpha-2)$ and put
$$
q=\frac{(\alpha-\eta)d}{2}>d.
$$
Recall that $R=Ma_\Lambda$ and $S_\Lambda=\Lambda a_\Lambda^d$.
We first estimate the two exterior contributions and then compare
the Riesz mean on $B_R$ with the sum over the complete cubes.

Let $r=a_\Lambda s$ with $s\geq1$, and take $\Lambda$ sufficiently
large that $a_\Lambda\geq e$. Since
$$
a_\Lambda\leq e+r\leq 2r=2a_\Lambda s,
$$
we have
$$
1\leq \frac{\log(e+r)}{\log a_\Lambda}
\le \frac{\log(2a_\Lambda s)}{\log a_\Lambda}
= 1+\frac{\log(2s)}{\log a_\Lambda}.
$$
If $\beta<0$, by $r>a_\Lambda$ and $\left(\frac{\log(e+r)}{\log a_\Lambda}\right)\ge 1$ we have
$$
\left(\frac{\log(e+r)}{\log a_\Lambda}\right)^\beta\le 1\le s^\eta,
\qquad\text{where}\ s\geq1,\ \eta\in(0,\alpha-2).
$$
If $\beta>0$, then $\log a_\Lambda\geq1$ gives
$$
\left(\frac{\log(e+r)}{\log a_\Lambda}\right)^\beta
\leq (1+\log(2s))^\beta
\leq C_{\beta,\eta}s^\eta,
\qquad s\geq1.
$$
The last inequality follows because $\frac{(1+\log(2s))^\beta}{s^\eta}$ is continuous for $s\geq1$ and tends to zero as $s\to\infty$. Thus, $\forall \beta \in \mathbb{R}$,
\begin{equation}\label{log}
\left(\frac{\log(e+r)}{\log a_\Lambda}\right)^\beta\leq C_{\beta,\eta}s^\eta,
\qquad s\geq1,\ \eta\in(0,\alpha-2).
\end{equation}
Using \eqref{log}, the definition of $w$ and the relation $\Lambda a_\Lambda^{-\alpha}(\log a_\Lambda)^\beta\to1$ as $\Lambda\rightarrow\infty$, we obtain
$$
\begin{aligned}
\Lambda w(r)
&=
\bigl[\Lambda a_\Lambda^{-\alpha}(\log a_\Lambda)^\beta\bigr]
\left(\frac{1+r}{a_\Lambda}\right)^{-\alpha}
\left(\frac{\log(e+r)}{\log a_\Lambda}\right)^\beta\\
&\leq C_{\alpha,\beta,\eta}s^{-\alpha+\eta},
\qquad s=\frac{r}{a_\Lambda}\geq1,
\end{aligned}
$$
where we used $\frac{(1+r)}{a_\Lambda}\geq s$. Consequently,
\begin{equation}\label{ww1}
(\Lambda W(x))^{\frac{d}{2}}
\leq C_{d,\alpha,\beta,\eta}
\left(\frac{|x|}{a_\Lambda}\right)^{-q},
\qquad |x|\geq a_\Lambda,
\end{equation}
where $q=\frac{(\alpha-\eta)d}{2}>d$.

For $M\geq1$ and $m\in \mathbb{N}$, define
$$
\mathcal A_m
=
\{x\in\mathbb Z^d:2^mMa_\Lambda\leq |x|<2^{m+1}Ma_\Lambda\}.
$$
The lattice point bound for balls gives
\begin{equation}\label{A1}
|\mathcal A_m|\leq C_d(2^mMa_\Lambda)^d.
\end{equation}
On $\mathcal A_m$, we also have $\frac{|x|}{a_\Lambda}\geq2^mM$. Since
$\frac{\Lambda^{1+\frac{d}{2}}}{S_\Lambda}=\frac{\Lambda^{\frac{d}{2}}}{a_\Lambda^d}$,
it follows \eqref{ww1} and \eqref{A1} that
\begin{equation}\label{W1}
\begin{aligned}
\frac{\Lambda^{1+\frac{d}{2}}}{S_\Lambda}
\sum_{|x|\geq Ma_\Lambda}W(x)^{\frac{d}{2}}
&=
\frac{1}{a_\Lambda^d}
\sum_{m=0}^{\infty}\sum_{x\in\mathcal A_m}
(\Lambda W(x))^{\frac{d}{2}}\\
&\leq
\frac{C_{d,\alpha,\beta,\eta}}{a_\Lambda^d}
\sum_{m=0}^{\infty}
(2^mMa_\Lambda)^d(2^mM)^{-q}\\
&=
C_{d,\alpha,\beta,\eta}M^{d-q}
\sum_{m=0}^{\infty}2^{m(d-q)}\\
&\leq C_{d,\alpha,\beta,\eta}M^{d-q}.
\end{aligned}
\end{equation}
The series converges because $q>d$, and the constant in the last
bound is independent of $M\geq1$ and $\Lambda$.

Dividing \eqref{R3} by $S_\Lambda$ and using \eqref{W1} gives
\begin{equation}\label{weyl-tail-error}
0\leq
\frac{\mathcal{R}(\Lambda)-\mathcal{R}_{B_R}(\Lambda)}{S_\Lambda}
\leq
C_d\frac{M^{d-1}}{a_\Lambda}
+
C_{d,\alpha,\beta,\eta}M^{d-q},
\end{equation}
because $\frac{\Lambda R^{d-1}}{S_\Lambda}=\frac{M^{d-1}}{a_\Lambda}$, where $q=\frac{(\alpha-\eta)d}{2}>d$.

The same exterior weight sum controls the tail of the lattice approximation. Indeed, Lemma~\ref{num:symbol} gives
$\rho(t)\leq C_dt^{1+\frac{d}{2}}$, so
\begin{equation}\label{W2}
0\leq W(x)^{-1}\rho(\Lambda W(x))
\leq C_d\Lambda^{1+\frac{d}{2}}W(x)^{\frac{d}{2}}.
\end{equation}
Combining \eqref{W1} and \eqref{W2}, we obtain
$$
\begin{aligned}
0
&\leq
\sum_{x\in\mathbb Z^d}W(x)^{-1}\rho(\Lambda W(x))
-
\sum_{x\in B_R}W(x)^{-1}\rho(\Lambda W(x))\\
&=
\sum_{|x|\geq R}W(x)^{-1}\rho(\Lambda W(x))\\
&\leq
C_d\Lambda^{1+\frac{d}{2}}\sum_{|x|\geq R}W(x)^{\frac{d}{2}}\\
&\leq C_{d,\alpha,\beta,\eta}M^{d-q}S_\Lambda.
\end{aligned}
$$

Applying the rank estimate from Lemma~\ref{num:block-rank} after multiplication
on both sides by $W_{B_R}^{-\frac{1}{2}}$ gives the same rank bound for
the weighted matrices. By Proposition~\ref{prop:weyl-rank}(iii),
and including the contribution $\Lambda|D_{R,L}|$ of the zero
block, we obtain
$$
\left|\mathcal{R}_{B_R}(\Lambda)
-\sum_{j\in J_{R,L}}\mathcal{R}_j^W(\Lambda)\right|
\leq C_d\Lambda\left(\frac{R^d}{L}+R^{d-1}L\right)
=o_M(S_\Lambda).
$$
Lemma~\ref{num:frozen} permits replacement of the variable weights by the
frozen weights. Lemma~\ref{num:periodic} and Lemma~\ref{num:lattice} then yield
$$
\mathcal{R}_{B_R}(\Lambda)
=
\sum_{x\in B_R}W(x)^{-1}\rho(\Lambda W(x))
+o_M(S_\Lambda),
$$
where $\frac{o_M(S_\Lambda)}{S_\Lambda}\to0$ as $\Lambda\to\infty$
with $M$ fixed.

By the triangle inequality and $\mathcal{R}(\Lambda)\geq \mathcal{R}_{B_R}(\Lambda)$,
$$
\begin{aligned}
&\left|\mathcal{R}(\Lambda)
-\sum_{x\in\mathbb Z^d}W(x)^{-1}\rho(\Lambda W(x))\right|\\
&\quad\leq \mathcal{R}(\Lambda)-\mathcal{R}_{B_R}(\Lambda)\\
&\qquad+
\left|\mathcal{R}_{B_R}(\Lambda)
-\sum_{x\in B_R}W(x)^{-1}\rho(\Lambda W(x))\right|\\
&\qquad+
\sum_{|x|\geq R}W(x)^{-1}\rho(\Lambda W(x)).
\end{aligned}
$$
Dividing by $S_\Lambda$ and using \eqref{weyl-tail-error}
for the first term and \eqref{W1} and \eqref{W2} for the last
term gives
$$
\begin{aligned}
\frac{1}{S_\Lambda}
\left|\mathcal{R}(\Lambda)
-\sum_{x\in\mathbb Z^d}W(x)^{-1}\rho(\Lambda W(x))\right|
\leq
C_d\frac{M^{d-1}}{a_\Lambda}
+C_{d,\alpha,\beta,\eta}M^{d-q}
+o_M(1),
\end{aligned}
$$
where $q=\frac{(\alpha-\eta)d}{2}>d$. Taking the limit superior as $\Lambda\to\infty$ with $M$ fixed
gives the upper bound $C_{d,\alpha,\beta,\eta}M^{d-q}$.
Since $q>d$, letting $M\to\infty$ proves the lemma.
\end{proof}
Next we finish the proof of Proposition~\ref{prop:riesz-asymptotic}.\\
\\
$\mathbf{The\ proof\ of\ the\ Proposition\ \text{\ref{prop:riesz-asymptotic}}}$ Set $h=a_\Lambda^{-1}$, so that $h\to0$ as $\Lambda\to\infty$.
For $y\in\mathbb R^d$, define
\begin{equation}\label{th}
t_h(y)=
\begin{cases}
\Lambda w(\frac{|y|}{h}), & y\ne0,\\
\Lambda W(0), & y=0,
\end{cases}
\qquad
F_h(y)=\frac{\rho(t_h(y))}{t_h(y)},
\end{equation}
where $w(r)=(1+r)^{-\alpha}[\log(e+r)]^\beta$.
Then $t_h(hx)=\Lambda W(x)$ for every $x\in\mathbb Z^d$, including the origin. Since $S_\Lambda=\Lambda a_\Lambda^d$, we have
$$
\frac{1}{S_\Lambda}
\sum_{x\in\mathbb Z^d} W(x)^{-1}\rho(\Lambda W(x))
=
h^d\sum_{x\in\mathbb Z^d}F_h(hx).
$$

For $y\ne0$, define
$$
F(y)=\frac{\rho(|y|^{-\alpha})}{|y|^{-\alpha}}
=|y|^\alpha\rho(|y|^{-\alpha}).
$$
Set $F(0)=0$. Since \eqref{th}, $0\leq\rho(t)\leq t$ and $\rho(t)\leq C_dt^{1+\frac{d}{2}}$, we have
$$
0\leq F_h(y)\leq1,
\qquad
0\leq F(y)\leq\min\{1,C_d|y|^{-\frac{\alpha d}{2}}\}
\quad (y\ne0).
$$
The bound on $F$ gives integrability because $\frac{\alpha d}{2}>d$.
Moreover, $\rho(t)>0$ for $t>4d$, since $\sigma(\theta)\leq4d$.
Hence $F(y)>0$ for $0<|y|<(4d)^{-\frac{1}{\alpha}}$, and
$$
0<K_{d,\alpha}=\int_{\mathbb R^d}F(y)\,dy<\infty.
$$

Fix $0<\varepsilon<1<T$. We divide the lattice sum and the
integral into three regions.

\medskip
\noindent\textit{Region 1: $|y|\leq\varepsilon$.}
The bounds $F_h\leq1$ and $F\leq1$ give
$$
h^d\sum_{|hx|\leq\varepsilon}F_h(hx)
\leq C_d(h+\varepsilon)^d,
\qquad
\int_{|y|\leq\varepsilon}F(y)\,dy
\leq C_d\varepsilon^d.
$$
The first estimate follows from the lattice-point bound
$\#\{x\in\mathbb Z^d:|x|\leq\frac{\varepsilon}{h}\}
\leq C_d(1+\frac{\varepsilon}{h})^d$.

\medskip
\noindent\textit{Region 2: $\varepsilon<|y|\leq T$.}
Since $h=a_\Lambda^{-1}$ and $t_h(hx)=\Lambda W(x)$, we can write
$$
t_h(y)=
\bigl[\Lambda a_\Lambda^{-\alpha}(\log a_\Lambda)^\beta\bigr]
(h+|y|)^{-\alpha}
\left(\frac{\log(e+a_\Lambda|y|)}{\log a_\Lambda}\right)^\beta.
$$
The first factor tends to $1$ as $\Lambda\rightarrow\infty$. The second factor tends uniformly to $|y|^{-\alpha}$ on the annulus as $\Lambda\rightarrow\infty$. For the last factor, observe that
$$
\frac{\log(e+a_\Lambda|y|)}{\log a_\Lambda}
=
1+\frac{\log(|y|+eh)}{\log a_\Lambda}\rightarrow1\quad\text{as }\Lambda\rightarrow\infty,
$$
where $|y|+eh$ is bounded. Consequently, $t_h(y)\to |y|^{-\alpha}$ uniformly on $\varepsilon<|y|\leq T$ as $\Lambda\rightarrow\infty$.

For all sufficiently small $h$, both $t_h(y)$ and
$|y|^{-\alpha}$ lie in a fixed compact subinterval of $(0,\infty)$
on this annulus. The uniform continuity of $t\mapsto\frac{\rho(t)}{t}$
on that interval therefore gives
$$
\sup_{\varepsilon<|y|\leq T}|F_h(y)-F(y)|\rightarrow0\quad\text{as }\Lambda\rightarrow\infty.
$$
Since $h^d\#\{x:|hx|\leq T\}\leq C_d(T+h)^d$, it follows that
$$
h^d\sum_{\varepsilon<|hx|\leq T}\bigl|F_h(hx)-F(hx)\bigr|\rightarrow0\quad\text{as }\Lambda\rightarrow\infty.
$$
The function $F$ is continuous on the closed annulus
$\varepsilon\leq|y|\leq T$, and its boundary spheres have
Lebesgue measure zero. Thus the Riemann sums give
$$
h^d\sum_{\varepsilon<|hx|\leq T}F_h(hx)
\rightarrow
\int_{\varepsilon<|y|\leq T}F(y)\,dy\quad\text{as }\Lambda\rightarrow\infty.
$$

\medskip
\noindent\textit{Region 3: $|y|>T$.}
Fix $\eta\in(0,\alpha-2)$ and put $q=\frac{(\alpha-\eta)d}{2}>d$.
Using \eqref{W1} and \eqref{W2} with $M=T$, we obtain
$$
\begin{aligned}
h^d\sum_{|hx|>T}F_h(hx)
&=
\frac{1}{S_\Lambda}
\sum_{|x|>Ta_\Lambda}W(x)^{-1}\rho(\Lambda W(x))\\
&\leq
\frac{C_d\Lambda^{1+\frac{d}{2}}}{S_\Lambda}
\sum_{|x|\geq Ta_\Lambda}W(x)^{\frac{d}{2}}\\
&\leq C_{d,\alpha,\beta,\eta}T^{d-q}.
\end{aligned}
$$
Since $0\leq F(y)\leq\min\{1,C_d|y|^{-\frac{\alpha d}{2}}\}$, we have
$$
\int_{|y|>T}F(y)\,dy
\leq \widetilde C_{d,\alpha}T^{d-\frac{\alpha d}{2}}.
$$

Combining the three regions and taking the limit superior as
$\Lambda\to\infty$ with $\varepsilon$ and $T$ fixed, we obtain
$$
\limsup_{\Lambda\to\infty}
\left|
h^d\sum_{x\in\mathbb Z^d}F_h(hx)-K_{d,\alpha}
\right|
\leq
C_d\varepsilon^d
+C_{d,\alpha,\beta,\eta}T^{d-q}
+\widetilde C_{d,\alpha}T^{d-\frac{\alpha d}{2}}.
$$
Since $q>d$ and $\frac{\alpha d}{2}>d$, letting $\varepsilon\rightarrow0$ and then $T\to\infty$ yields
$$
\frac{1}{S_\Lambda}
\sum_{x\in\mathbb Z^d}W(x)^{-1}\rho(\Lambda W(x))
\rightarrow K_{d,\alpha}\quad\text{as }\Lambda\rightarrow\infty.
$$
By Lemma~\ref{num:full}, the difference between $\mathcal{R}(\Lambda)$ and the lattice sum is $o(S_\Lambda)$. Therefore
$$
\mathcal{R}(\Lambda)\sim K_{d,\alpha}S_\Lambda
=K_{d,\alpha}\Lambda a_\Lambda^d\quad\text{as }\Lambda\rightarrow\infty.
$$
\begin{flushright}
  $\Box$
\end{flushright}

\subsection{Proof of Theorem 1.3}\label{num:proof-weyl}
Set $p=1+\frac{d}{\alpha}$ and $\widehat C_{d,\alpha}=pK_{d,\alpha}$, so that $\widehat C_{d,\alpha}=C_{d,\alpha}^{-\frac{d}{\alpha}}$ with $C_{d,\alpha}$ as in \eqref{eq Cdef}. Since $\mathcal{R}(\Lambda)=\int_0^\Lambda \mathcal{N}(t)\,dt$ and $\mathcal{N}(\Lambda)$ is nondecreasing, for $0<s<\Lambda$, by the monotonicity of $\mathcal{N}$ we have
$$
\mathcal{R}(\Lambda)-\mathcal{R}(\Lambda-s)=\int_{\Lambda-s}^{\Lambda}\mathcal{N}(t)\,dt\leq s\mathcal{N}(\Lambda)
$$
and
$$
\mathcal{R}(\Lambda+s)-\mathcal{R}(\Lambda)=\int_{\Lambda}^{\Lambda+s}\mathcal{N}(t)\,dt\geq s\mathcal{N}(\Lambda).
$$
Consequently,
\begin{equation}\label{N1}
\frac{\mathcal{R}(\Lambda)-\mathcal{R}(\Lambda-s)}{s}\leq \mathcal{N}(\Lambda)\leq\frac{\mathcal{R}(\Lambda+s)-\mathcal{R}(\Lambda)}{s}.
\end{equation}
Fix $\varepsilon\in(0,1)$, and take $s=\varepsilon\Lambda$. Then
$$
\frac{\mathcal{R}(\Lambda)-\mathcal{R}((1-\varepsilon)\Lambda)}
{\varepsilon\Lambda}
\leq \mathcal{N}(\Lambda)
\leq
\frac{\mathcal{R}((1+\varepsilon)\Lambda)-\mathcal{R}(\Lambda)}
{\varepsilon\Lambda}.
$$
For each fixed $c>0$, the definition of $a_\Lambda$ gives
$\frac{a_{c\Lambda}}{a_\Lambda}\to c^{\frac{1}{\alpha}}$ as $\Lambda\rightarrow\infty$. By Proposition~\ref{prop:riesz-asymptotic},
$$
\frac{\mathcal{R}(c\Lambda)}{\Lambda a_\Lambda^d}
\rightarrow K_{d,\alpha}c^p\quad\text{as }\Lambda\rightarrow\infty.
$$
Dividing the difference-quotient bounds by $a_\Lambda^d$ and
letting $\Lambda\to\infty$, we obtain
$$
K_{d,\alpha}\frac{1-(1-\varepsilon)^p}{\varepsilon}
\leq \liminf_{\Lambda\to\infty}\frac{\mathcal{N}(\Lambda)}{a_\Lambda^d}
\leq \limsup_{\Lambda\to\infty}\frac{\mathcal{N}(\Lambda)}{a_\Lambda^d}
\leq K_{d,\alpha}\frac{(1+\varepsilon)^p-1}{\varepsilon}.
$$
Besides, since
$$
\frac{1-(1-\varepsilon)^p}{\varepsilon}
\rightarrow p,
\qquad
\frac{(1+\varepsilon)^p-1}{\varepsilon}
\rightarrow p\quad\text{as }\varepsilon\rightarrow0,
$$
where $p=1+\frac{d}{\alpha}$. Thus let $\varepsilon\rightarrow0$ we have
\begin{equation}\label{NN1}
\mathcal{N}(\Lambda)\sim \widehat C_{d,\alpha}a_\Lambda^d
=
(1+\frac{d}{\alpha})K_{d,\alpha}a_\Lambda^d\quad\text{as }\Lambda\rightarrow\infty.
\end{equation}
We now use $W$ for the general weight in Theorem~1.3.
Since $\frac{W(x)}{(1+|x|)^{-2}}\to0$ as $|x|\to\infty$, the same
compactness argument applies to $W$. Write $\mathcal{N}_W$ and
$\mathcal{N}_{W_0}$ for the counting functions of the two weighted
problems, with eigenvalues counted with multiplicity.

For each $\varepsilon\in(0,1)$, choose $R_\varepsilon>1$ such
that, for $|x|>R_\varepsilon$,
$$
(1-\varepsilon)W_0(x)\leq W(x)\leq(1+\varepsilon)W_0(x).
$$
Such a radius exists because $\frac{W(x)}{W_0(x)}\to1$ as
$|x|\to\infty$. Set
$$
B_{R_\varepsilon}(0)=\{x\in\mathbb{Z}^d:|x|\leq R_\varepsilon\}.
$$
For a weight $W$ and $W_0$, write
$$
\|u\|_W^2=\sum_{x\in\mathbb{Z}^d}W(x)|u(x)|^2\ \text{and}\ \|u\|_{W_0}^2=\sum_{x\in\mathbb{Z}^d}W_0(x)|u(x)|^2.
$$
By the min--max principle, $\mathcal N_S(t)$ is the largest dimension of a subspace on which every function $u$ satisfies $\|u\|_1^2\leq t\|u\|_S^2$, where $S=W$ or $W_0$.

We first prove the lower bound for $\mathcal N_W(\Lambda)$.
Let $m=\mathcal N_{W_0}((1-\varepsilon)\Lambda)$, and let
$\xi_1,\ldots,\xi_m$ be the corresponding eigenfunctions,
orthonormal with respect to $W_0$. Denote their eigenvalues
by $\lambda_1,\ldots,\lambda_m$. For every linear combination
$u=\sum_{k=1}^m c_k\xi_k$, we have
$$
\|u\|_1^2
=\sum_{k=1}^m\lambda_k|c_k|^2
\leq (1-\varepsilon)\Lambda\sum_{k=1}^m|c_k|^2
=(1-\varepsilon)\Lambda\|u\|_{W_0}^2.
$$

We now choose the coefficients so that $u(x)=0$ at every $x\in B_{R_\varepsilon}(0)$. At each such point, this condition is the linear equation
$$
\sum_{k=1}^m c_k\xi_k(x)=0.
$$
There are $m$ coefficients and $|B_{R_\varepsilon}(0)|$ equations, so the space of coefficient vectors satisfying all these equations has dimension at least $m-|B_{R_\varepsilon}(0)|$. Since the $\xi_k$ are linearly independent,
the corresponding linear combinations form a subspace of the
same dimension.

For every such linear combination, all terms with $x\in B_{R_\varepsilon}(0)$ vanish in both weighted sums. Hence the comparison of the weights outside $B_{R_\varepsilon}(0)$ gives
$$
\begin{aligned}
\|u\|_W^2
&=\sum_{x\notin B_{R_\varepsilon}(0)}W(x)|u(x)|^2\\
&\geq (1-\varepsilon)
      \sum_{x\notin B_{R_\varepsilon}(0)}W_0(x)|u(x)|^2\\
&=(1-\varepsilon)\|u\|_{W_0}^2.
\end{aligned}
$$
These functions are still linear combinations of
$\xi_1,\ldots,\xi_m$, so the preceding energy estimate yields
$$
\|u\|_1^2\leq(1-\varepsilon)\Lambda\|u\|_{W_0}^2
\leq\Lambda\|u\|_W^2.
$$
By the min--max principle, it follows that
$$
\mathcal N_W(\Lambda)
\geq \mathcal N_{W_0}((1-\varepsilon)\Lambda)-|B_{R_\varepsilon}(0)|.
$$
For the upper bound, by the same method we have
$$
\mathcal N_{W_0}((1+\varepsilon)\Lambda)
\geq \mathcal N_W(\Lambda)-|B_{R_\varepsilon}(0)|.
$$
Combining the lower and upper bounds, we obtain
$$
\mathcal N_{W_0}((1-\varepsilon)\Lambda)-|B_{R_\varepsilon}(0)|
\leq \mathcal N_W(\Lambda)
\leq \mathcal N_{W_0}((1+\varepsilon)\Lambda)+|B_{R_\varepsilon}(0)|.
$$

For fixed $\varepsilon$, the set $B_{R_\varepsilon}(0)$ is finite and independent of $\Lambda$, so $\frac{|B_{R_\varepsilon}(0)|}{a_\Lambda^d}\to0$. Dividing by $a_\Lambda^d$ and using \eqref{NN1} for the model weight, together with $\frac{a_{c\Lambda}}{a_\Lambda}\to c^{\frac{1}{\alpha}}$ for each fixed $c>0$, we obtain
$$
\begin{aligned}
\widehat C_{d,\alpha}(1-\varepsilon)^{\frac{d}{\alpha}}
&\leq \liminf_{\Lambda\to\infty}
       \frac{\mathcal N_W(\Lambda)}{a_\Lambda^d}\\
&\leq \limsup_{\Lambda\to\infty}
       \frac{\mathcal N_W(\Lambda)}{a_\Lambda^d}
\leq \widehat C_{d,\alpha}(1+\varepsilon)^{\frac{d}{\alpha}}.
\end{aligned}
$$
Letting $\varepsilon\rightarrow0$ gives
$$
\mathcal N_W(\Lambda)\sim \widehat C_{d,\alpha}a_\Lambda^d
\quad\text{as }\Lambda\to\infty.
$$
In the remainder of the proof, $\mathcal N$ denotes
$\mathcal N_W$.
To recover the eigenvalue asymptotics from the counting function, set
\begin{equation}\label{b1}
b_k
=
d^\beta
\left(\frac{k}{\widehat C_{d,\alpha}}\right)^{\frac{\alpha}{d}}
(\log k)^{-\beta}\quad k\geq 2.
\end{equation}
Taking logarithms gives
$$
\log b_k
=
\frac{\alpha}{d}\log k
-\beta\log\log k
+\beta\log d
-\frac{\alpha}{d}\log \widehat C_{d,\alpha}.
$$
Since $\log\log k=o(\log k)$, we have
\begin{equation}\label{loglog}
\log b_k\sim(\frac{\alpha}{d})\log k\quad \text{as}\ k\rightarrow\infty.
\end{equation} In particular, $b_k\to\infty$ as $k\rightarrow\infty$.
Applying the counting asymptotic \eqref{NN1} at $\Lambda=b_k$, \eqref{b1} and \eqref{loglog} we obtain
\begin{equation}\label{FF1}
\begin{aligned}
\frac{\mathcal{N}(b_k)}{k}
&\sim
\widehat C_{d,\alpha}\alpha^{-\frac{\beta d}{\alpha}}
\frac{b_k^{\frac{d}{\alpha}}}{k}(\log b_k)^{\frac{\beta d}{\alpha}}\\
&=
\left(\frac{d}{\alpha}\frac{\log b_k}{\log k}\right)^{\frac{\beta d}{\alpha}}
\rightarrow 1\quad\text{as }k\rightarrow\infty.
\end{aligned}
\end{equation}
For every fixed $c>0$, we have $\frac{\log(cb_k)}{\log b_k}\to1$ as $k\rightarrow\infty$. A second application of the counting asymptotic \eqref{NN1} gives
\begin{equation}\label{FF2}
\frac{\mathcal{N}(cb_k)}{\mathcal{N}(b_k)}
\sim
c^{\frac{d}{\alpha}}
\left(\frac{\log(cb_k)}{\log b_k}\right)^{\frac{\beta d}{\alpha}}
\rightarrow c^{\frac{d}{\alpha}}\quad\text{as }k\rightarrow\infty.
\end{equation}
Consequently, by \eqref{FF1} and \eqref{FF2} we have
$$
\frac{\mathcal{N}(cb_k)}{k}\to c^{\frac{d}{\alpha}}\quad \text{as}\ k\rightarrow\infty.
$$
Fix $\varepsilon\in(0,1)$. Since $(1-\varepsilon)^{\frac{d}{\alpha}}<1<(1+\varepsilon)^{\frac{d}{\alpha}}$, for all sufficiently large $k$, take $c$ equal to $1-\varepsilon$ and $1+\varepsilon$ respectively, we have
$$
\mathcal{N}((1-\varepsilon)b_k)<k<\mathcal{N}((1+\varepsilon)b_k).
$$
It follows that
$$
(1-\varepsilon)b_k\leq\lambda_k\leq(1+\varepsilon)b_k.
$$
Letting $\varepsilon\rightarrow0$ and by \eqref{b1} we have
$$
\lambda_k\sim
d^\beta C_{d,\alpha}
k^{\frac{\alpha}{d}}
(\log k)^{-\beta}
\quad\text{as }k\to\infty.
$$
\hfill$\Box$

  \bigskip
  \medskip

 \small
\noindent {\bf  Conflicts of interest:} The authors declare that they have no conflicts of interest regarding this work.
\medskip

\noindent {\bf  Data availability:}  This paper has no associated data.\medskip

\noindent{\bf Acknowledgements:}
 H. Chen is supported by  the National Natural Science Foundation of China, No.  12361043.
\end{document}